\documentclass[12pt]{iopart}

\usepackage{graphicx,cite,ulem}
\usepackage{subfig}
\usepackage{multirow}
\usepackage{xcolor}
\usepackage{float}
\usepackage{placeins}

\renewcommand{\em}{\it}
	\newcommand{\ignore}[1]{}
\newcommand{\eqref}{\eref}

\begin{document}

\title{Pulsating dynamics of slow-fast population models with delay}
\author{Pavel Kravetc$^1$, Dmitrii Rachinskii$^1$, Andrei Vladimirov$^{2,3}$}
\address{$^1$ Department of Mathematical Sciences, The University of Texas at Dallas, USA}
\address{$^2$ Weierstrass Institute, Mohrenstr. 39, 10117 Berlin, Germany}
\address{$^3$ Lobachevsky State University of Nizhny Novgorod, Nizhny Novgorod, Russia}

\ead{pxk142530@utdallas.edu, dmitry.rachinskiy@utdallas.edu and vladimir@wias-berlin.de}

\begin{abstract}
We discuss a bifurcation scenario which creates periodic pulsating solutions in slow-fast delayed systems through a cascade
of almost simultaneous Hopf bifurcations.
This scenario has been previously associated with formation of pulses in a delayed model of mode-locked semiconductor lasers.
In this work, through a case study of several examples, we establish that a cascade of Hopf bifurcations
can produce periodic pulses, with a period close to the delay time, in population dynamics models and explore the conditions
that ensure the realization of this scenario near a transcritical bifurcation threshold.
We derive asymptotic approximations for the pulsating solution and consider scaling of the solution and its period
with the small parameter that measures the ratio of the time scales. The role of competition for the realization
of the bifurcation scenario is highlighted.
\\


\end{abstract}


\maketitle

\section{Introduction}\label{intro}

Mode locking of lasers \cite{Haus00} is a powerful technique used to produce a periodic sequences of short optical pulses at high repetition rates, which are suitable for various applications, including material processing, medical imaging, telecommunications \cite{Ippen,Kaiser}, optical sampling, microwave photonics, optical division multiplexing \cite{Delfyett}, and two-photon imaging \cite{Kuramoto}.
The optical spectrum of a mode-locked laser consists of a set of equally spaced narrow lines corresponding to the longitudinal cavity modes characterized by certain fixed phase relationships between them. Achieving such phase relationships can be, at least qualitatively, viewed as a problem of synchronization of many nonlinear coupled oscillators with frequencies close to multiples of a fundamental frequency. There are two main methods to produce mode-locked optical pulses: active and passive mode-locking, and also a combination of the two techniques called hybrid mode-locking. In particular, a passively mode-locked laser is a self-oscillating system which does not require the use of an external radio frequency modulation. Passive mode-locking is commonly achieved by including a saturable absorber into the laser cavity. In the classical theory of a mode-locked laser 
due to Haus \cite{Haus}, a slow evolution of the shape of the optical pulse circulating in the cavity is described by a complex parabolic 
master equation of Ginzburg-Landau type. 
The solution describing a solitary pulse is explicit and has the hyperbolic secant profile.

The Haus master equation is derived under the assumption of small gain and loss per cavity round trip.
An alternative functional differential model, which is free from this approximation, has been obtained from the traveling wave model in the case of a ring geometry of the laser cavity in \cite{VladimirovTuraev}. Under further natural assumptions, such as the Lorentzian profile of the spectral filtering element, the functional differential model simplifies to a system of delay differential equations with a single delay corresponding to the cold cavity round trip time  \cite{VladimirovTuraev04,VladimirovTuraev,VladimirovTuraevKozyreff}. This delay differential model is suitable for describing mode-locking in a laser with large gain and losses, that is the situation typical of semiconductor laser devices\footnote{In the limit of small gain and losses per cavity round trip, one recovers
the Haus hyperbolic secant pulse shape in the delay differential model.}.
At the same time, the model is amenable to analytical and numerical bifurcation analysis \cite{VladimirovTuraevKozyreff,Rachinskii,VRW12,VRW15,VRW20}.

The delay differential model proposed in \cite{VladimirovTuraev04,VladimirovTuraevKozyreff,VladimirovTuraev} is a multi-rate system where several disparate time scales can be identified.
It has been shown that the transition from a relative equilibrium (continuous wave laser operation) to the mode-locked states with the increase of a bifurcation parameter
(the pump current) in this model is associated with a sequence of Hopf bifurcations on the relative equilibrium \cite{VladimirovTuraevKozyreff}.
The goal of this paper is to present an evidence that a cascade of almost simultaneous and almost resonant Hopf bifurcations is responsible for the formation of periodic pulsating dynamics with specific properties of a mode-locked regime for a broad class of slow-fast delayed differential systems. This evidence will be provided by a case study of predator-prey population dynamics models. An asymptotic analysis of the Hopf bifurcation cascade scenario will equip us with a tool for identifying and understanding conditions, which initialize the periodic pulsating dynamics, and for finding parameter values that support the scenario of a fast transition from a steady state to the periodic pulsation.

A delay is generally believed to be a destabilizing factor in population dynamics models \cite{AndersonMay}.
In particular, increasing delay can lead to oscillations where the system with zero or small delay exhibits a globally stable equilibrium. At the same time, pulsating solutions where time intervals of almost complete extinction of some species alternate with outbursts in their number are a typical feature of population dynamics models such as Lotka-Volterra systems, host-parasite models, susceptible-infective-recovered (SIR) epidemiological models and chemical kinetics models.

In this paper, we consider models involving populations of species which evolve on two different time scales. The models include an explicit delay time $T>0$ which can have different nature and, therefore, can appear in different terms of the equations \cite{Ruan}; the maturity delay is considered as the main example.
We are concerned with periodic dynamics presented by a limit cycle with the following properties:
\begin{itemize}
\item The period of the cycle is close to the delay time $T$;
\item The time trace of one component (which we call the $A$-compon\-ent) of the cycle is a sequence of identical short pulses, typically one pulse per period,
separated by intervals where the $A$-component is close to zero;
\item The oscillations are self-excited, {\em i.e.}, the cycle is either globally stable or has a large basin of attraction,
while the equilibrium with the zero $A$-component is unstable.
\end{itemize}
These properties will be formalized and quantified in terms of the parameter $\gamma\gg 1$ which measures the ratio of the slow and fast time scales of the population processes involved in the system. In particular, the period of the cycle is $T+O(1/\gamma)$, the duration of the pulse scales as $1/\gamma$, while the pulse amplitude is asymptotically proportional to $\gamma$, and the time average of each population tends to some limit value
as $\gamma$ increases.

We will look into a mechanism which generates a cycle with the above properties.
This mechanism will be associated with a cascade of Hopf bifurcations with almost commensurate frequencies
close to the multiples of $2\pi/T$, which all occur over a small interval of bifurcation parameter values. In the scenario we describe, this cascade is almost simultaneous with the transcritical bifurcation separating
the domain where the equilibrium with the zero $A$-component is stable from the domain
where it is unstable and coexists with the positive equilibrium. That is, the periodic pulsating dynamics
develops almost immediately after the bifurcation parameter crosses the threshold value (at the transcritical bifurcation point) via the cascade of Hopf bifurcations on the positive equilibrium.
This mechanism can be suspected whenever an observed/measured cyclic dynamics is represented by a
periodic sequence of pulses separated by intervals of almost complete extinction of some species and the period of this dynamics can be associated with some well defined delay in a system, which involves processes running on different time scales. The interplay between the delay and the slow-fast structure of the system
{is} at the heart of this scenario\footnote{The scenario we discuss is not related to switching between stable branches of
a critical manifold of a singularly perturbed system.}. We will obtain an approximation and important parameters of the pulses by adapting an approach proposed in \cite{New,VladimirovTuraev} for analysis
of optical systems to the population dynamics models.
We will also present an informal discussion of the mechanism creating the pulsating regime in
terms of the effect of increasing/decreas\-ing populations on each other. The role of the competition between the fast species will be considered in this context. The qualitative discussion should complement our attempt to understand and quantify the mechanism of mode-locking in terms of the associated bifurcation scenario.

The paper is organized as follows. In Section \ref{main}, we introduce a prototype population model and discuss the cascade of Hopf bifurcations, which we associate with the formation of a pulsating periodic solution.
In Section \ref{variations}, three variations of the prototype model are considered in order to demonstrate similar results for a wider class of population systems. In Section \ref{scaling}, we derive
 asymptotic approximations for the pulsating periodic solutions. Finally, in Section \ref{competing} the role of the competition for the realization of the bifurcation scenario is highlighted.

\section{Main prototype model}\label{main}

Our main prototype model has the form
\begin{eqnarray}
\gamma^{-1} A' &=& -A + \kappa G(t-T)A(t-T) - \mu Q A,\label{A}\\
\gamma^{-1} Q' &=& q_0-\beta Q -s A Q, \label{Q}\\
 G' &=& g_0-\alpha G - k A G \label{G}
\end{eqnarray}
where $T$ is the maturity delay of the species $A$, see \cite{Ruan}. Here $\gamma\gg1$, that is the species $A$ and $Q$ are assumed to be fast
(have much faster metabolism, higher reproductive rate etc.) compared to the species $G$. The species $A$ is a predator for the prey $G$.
The species $Q$ competes with $A$. All the parameters are positive\footnote{{The death rate of the species $A$ is scaled to 1. The number of parameters can be further reduced in a standard way by rescaling the phase variables and time.}}.

The species $Q$ plays an important role which will be clarified in further sections.
In particular, we will see that the system of the two equations \eqref{A} (with zero $Q$) and \eqref{G} does not demonstrate pulsating dynamics near the threshold.

The species $Q$ and $G$ are assumed to be recruited through constant immigration in Eqs. \eqref{A}--\eqref{G}.
In further sections, we will show that similar systems with different recruitment terms, including recruitment with constant birth rate,
show similar pulsating dynamics near the threshold. Also, delaying different terms seems to have little effect on solutions
in our examples; for instance, replacing the delayed term $G(t-T)$ by $G(t)$ in Eqs. \eqref{A}--\eqref{G}
preserves the periodic pulsating dynamics.

We will discuss nonnegative solutions only. Note that system \eqref{A}--\eqref{G}, as well as all the other systems considered in the paper, is positively invariant.

We associate the pulsating regime of system \eqref{A}--\eqref{G} near the point of the transcritical bifurcation of equilibria with
the Hopf bifurcations from the positive equilibrium. The recruitment rate $g_0$ of the prey $G$ will be used as the bifurcation parameter.

System \eqref{A}--\eqref{G} has an equilibrium with zero $A$,
\begin{equation}\label{equil1}
A_o=0, \qquad Q_o=\frac{q_0}{\beta},\qquad G_o=\frac{g_0}{\alpha},
\end{equation}
for all positive $g_0$, and a positive equilibrium
either for $g_0> g_0^*$ or for $g_0< g_0^*$, where the threshold value $g_0^*$ is defined by
\begin{equation}\label{thresho}
\frac{\kappa g_0^*}{\alpha} - \frac{\mu q_0}{\beta} = 1.
\end{equation}
{These} two equilibria collide in a transcritical bifurcation for $g_0=g_0^*$.
The positive equilibrium near the threshold is defined by the asymptotic formulas
\begin{equation}\label{equil2}
\fl A_*=\tilde a \delta +O(\delta^2), \qquad Q_*=\frac{q_0}{\beta} + \tilde q \delta+O(\delta^2),\qquad G_*=\frac{g_0^*}{\alpha}+\tilde g \delta+O(\delta^2)
\end{equation}
where $\delta=g_0-g_0^*$ and the coefficients of the first order correction are given by
\begin{equation*}
\tilde a=\frac 1 {\frac{k g^*_0}\alpha-\frac{\alpha \mu s q_0}{\kappa \beta^2}},\qquad \tilde q=\frac 1 {\frac {\alpha \mu} \kappa - \frac{k g_0^*\beta^2}{\alpha s q_0}}, \qquad \tilde g=\frac{\mu}{\kappa} \tilde q.
\end{equation*}
We will assume that
\begin{equation}\label{trans}
\frac{k g^*_0}{\alpha^2}>\frac{\mu s q_0}{\kappa \beta^2}.
\end{equation}
In this case, the positive equilibrium exists for $g_0>g_0^*$ and is stable near the threshold.
(If the opposite inequality holds, then the positive equilibrium exists for $g_0<g_0^*$ and is unstable near the threshold.)

The eigenvalues of the linearization of system \eqref{A}--\eqref{G} at the equilibrium \eqref{equil1} with zero $A$ are defined by
the relations $\lambda=-\gamma\beta<0$, $\lambda=-\alpha<0$ and
\begin{equation}\label{eig}
1+\frac{\lambda}{\gamma} = \frac{\kappa g_0}{\alpha} e^{-\lambda T} - \frac{\mu q_0}{\beta}.
\end{equation}
The solutions of \eqref{eig} satisfy ${\rm Re}\, \lambda<0$ in a left neighborhood of the threshold,
more precisely, for $g_0< g_0^* = \alpha(1+\mu q_0/\beta)/\kappa$.
Hence, the equilibrium \eqref{equil1} is stable below the threshold, {\em i.e.}, for $g_0<g_0^*$.
Consequently, the positive equilibrium \eqref{equil2} is stable in a small right neighborhood of the threshold,
{\em i.e.}, for small $\delta=g_0-g_0^*>0$.

\begin{figure}[ht!]
\begin{center}
\includegraphics*[width=0.5\columnwidth]{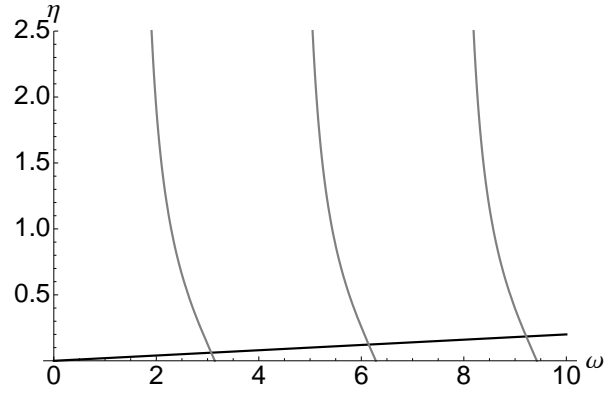}
\end{center}
\caption{\small{Solution of \eqref{transc}. The horizontal axis is $\omega$. Every second intersection of the straight line $\eta =\omega/\gamma $ and the function $\eta=-\kappa g_0^*\tan (\omega T)/\alpha$
satisfies the condition \eqref{transc1}. Here $\gamma=100$, $T=1$. \label{fig1}}}
\end{figure}

The equilibrium \eqref{equil1} undergoes a sequence of Hopf bifurcations in a small right neighborhood of
the threshold $g_0=g_0^*$ for large $\gamma$. To see this, first note that in the limit $\gamma = \infty$ the solutions of the characteristic equation \eqref{eig} have the form
\begin{equation*}
\lambda=i \omega_n, \qquad \omega_n=\frac{2\pi n}{T},\qquad n=1,2,\ldots,
\end{equation*}
{\em i.e.}, the equilibrium satisfies the necessary condition for
infinitely many simultaneous Hopf bifurcations at the threshold point $g_0=
g_0^*$. Moreover, these bifurcations are in resonance with each other as the frequencies
$\omega_n$ are all multiples of $2\pi/T$. For finite $\gamma$, setting $\lambda=i\omega$ in \eqref{eig}
in order to satisfy the Hopf bifurcation condition, and rearranging, we obtain the equations
\begin{eqnarray}\label{transc}
\frac{\omega}{\gamma} & = & - \frac{\kappa g_0^*}{\alpha} \tan {\omega T},
\\
\label{transc1}
\delta & = & g_0^*\left(\frac {1} {\cos{\omega T}} - 1 \right)>0
\end{eqnarray}
which define the frequency of the cycle and the bifurcation value of the parameter $g_0=g_0^*+\delta$
for each Hopf bifurcation from the equilibrium \eqref{equil1}.
Figure \ref{fig1} illustrates solutions of the transcendental equation \eqref{transc}.
For $\gamma\gg 1$, the solutions of Eqs. \eqref{transc}, \eqref{transc1} are approximated by the asymptotic formulas
\begin{equation}\label{bifi}
\omega_n = \frac{2\pi n}{T}\left(1-\frac {\alpha} {\kappa g_0^*\gamma T} + \frac {\alpha^2} {\left(\kappa g_0^* \gamma T\right)^2}\right) + O\left(\gamma^{-3}\right),
\end{equation}
\begin{equation}\label{bifi:delta}\delta_n = \frac {\alpha^2}{2 \kappa^2 g_0^*} \left(\frac{2 \pi n}{\gamma T} \right)^2 + O \left(\gamma^{-3} \right)
\end{equation}
with $n=1,2,\ldots$\ Hence, the $n$-th Hopf bifurcation after the threshold has a frequency close
to $2\pi n/T$ and $O(\sqrt{\gamma})$ Hopf bifurcations occur within the distance of order $1/\gamma$
from the threshold on the parameter $g_0$ axis.

\begin{figure}[ht!]
\begin{center}
\subfloat[\label{stabilZero}]{\includegraphics*[width=.5\textwidth]{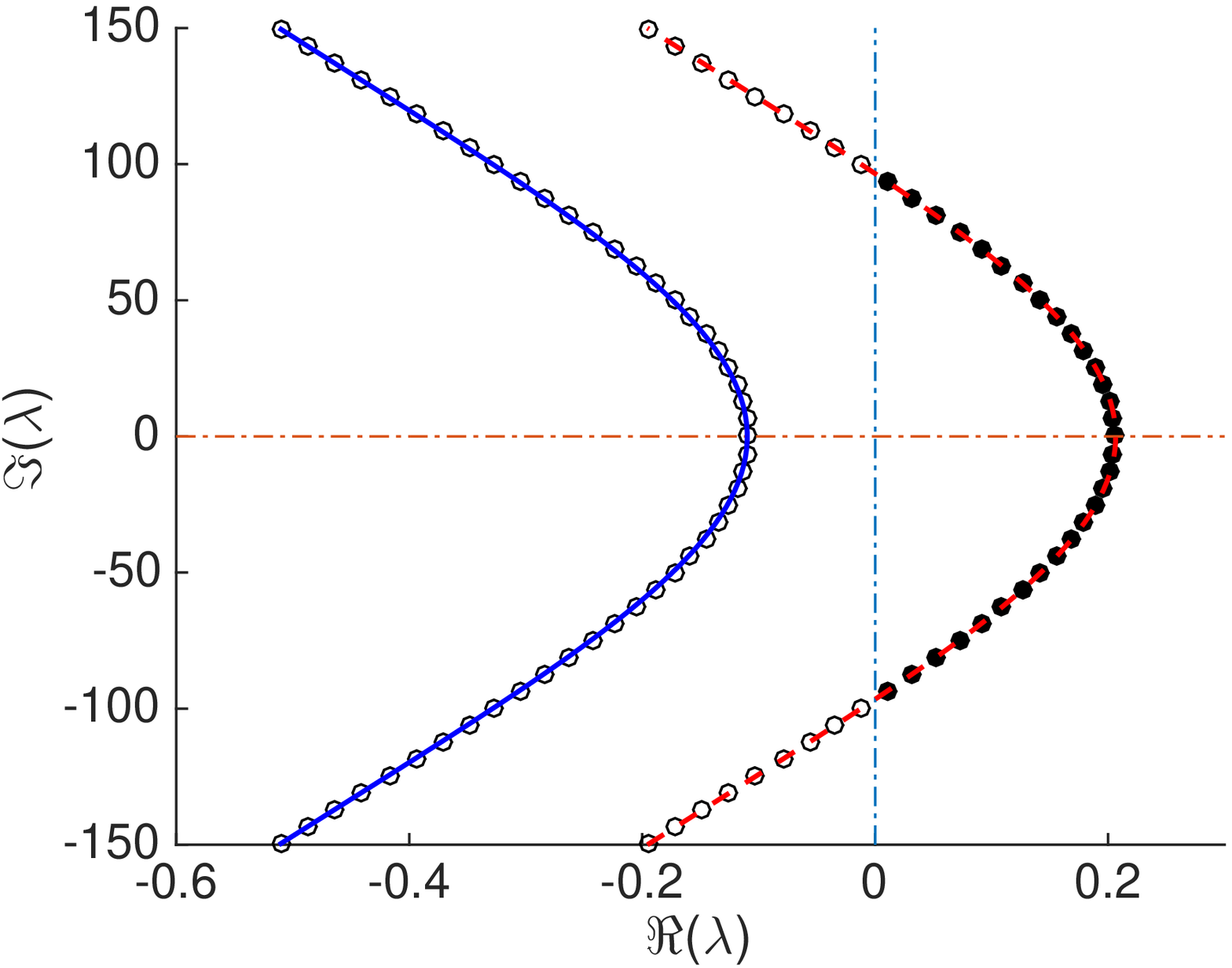}}
\hfill
\subfloat[\label{stabil50}]{\includegraphics*[width=.5\textwidth]{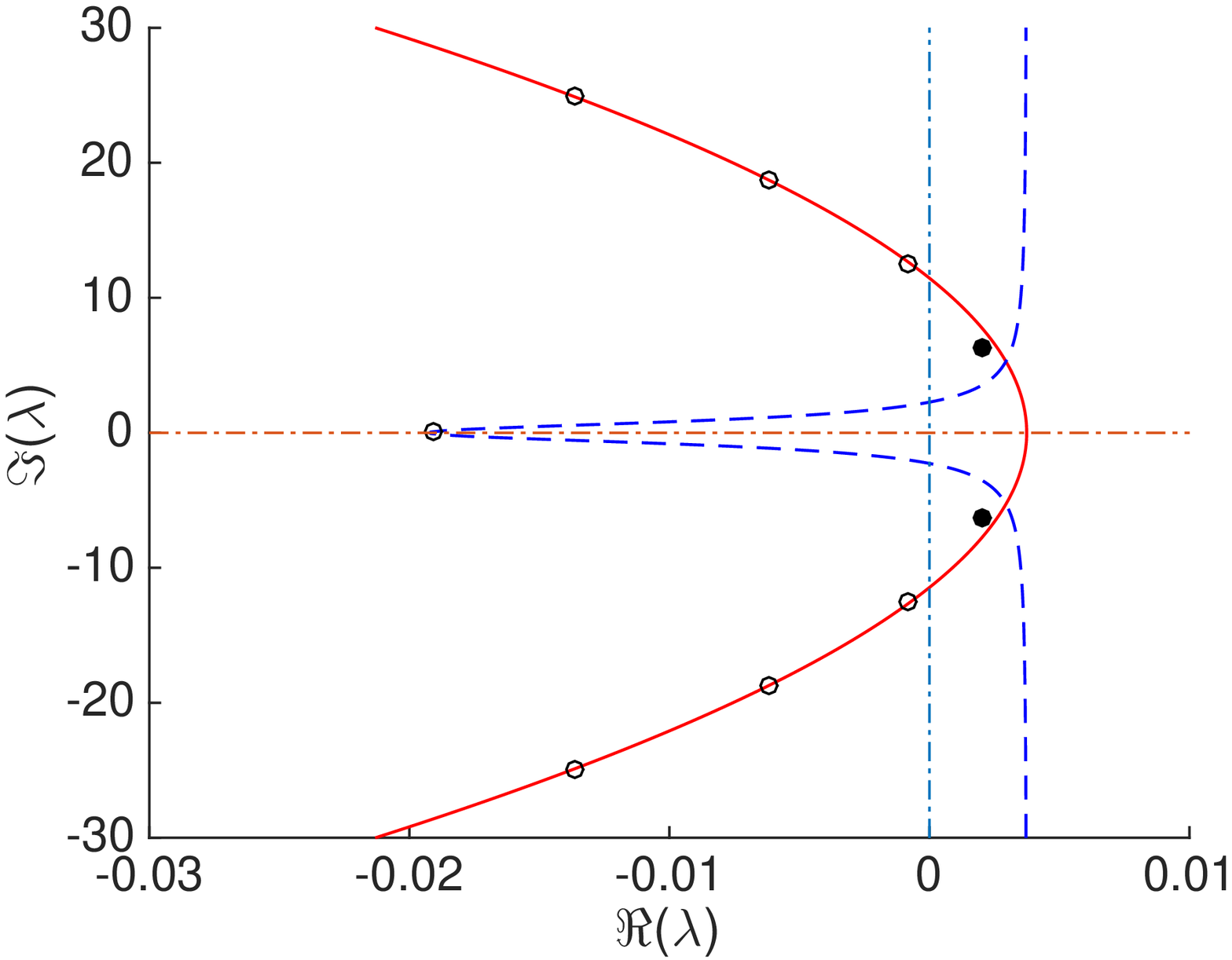}}
\end{center}
\caption{\small{Panel (a): Spectrum of the zero equilibrium for $g_0 = 2.6880$ and $g_0 = 3.6880$.
Numerical values of the eigenvalues are shown by circles; lines are obtained from \eqref{curve}.
Filled circles correspond to unstable eigenvalues.
Panel (b): Spectrum of the positive equilibrium of system \eqref{A}--\eqref{G} after the first Hopf bifurcation ($g_0 = 3.0269$), \textit{i.e.} exactly one pair of complex conjugate eigenvalues crossed the imaginary axis. Other parameters are the same as in Figure \ref{bifdiag2}.
Solid line defined by \eqref{weak} carries the weak spectrum; dashed line \eqref{strong} carries strongly stable spectrum.}}
\label{stabil}
\end{figure}

Following \cite{wy, yw}, the spectrum of the zero equilibrium defined by \eqref{eig} can be called {\sl weak} or {\sl pseudocontinuous} spectrum.
It is characterized by a specific scaling of the real and imaginary parts of the eigenvalues $\lambda = x + i \gamma \omega$
with $\gamma\gg 1$, where $x$ and $\omega$ are of order $1$. Using this scaling, we obtain from \eqref{eig} an approximate relationship between the real and imaginary parts
of the eigenvalues:
\begin{equation}\label{curve}
	x(\omega) = \frac{1}{2 T} \left(2\ln {\left(\frac{g_0}{g_0^*}\right)} - \ln{\left(1+\left(\frac{\alpha \omega}{ g_0^* \kappa}\right)^2\right)}\right)+O(\gamma^{-1}),
\end{equation}
which is dual to formulas \eqref{bifi}, \eqref{bifi:delta}.
The curve \eqref{curve} carrying the eigenvalues simply moves to the right with increasing $g_0$, see Figure \ref{stabilZero}.

\begin{figure}[ht]
\begin{center}
\subfloat[\label{bifdiag2}]{\includegraphics*[width=.5\textwidth]{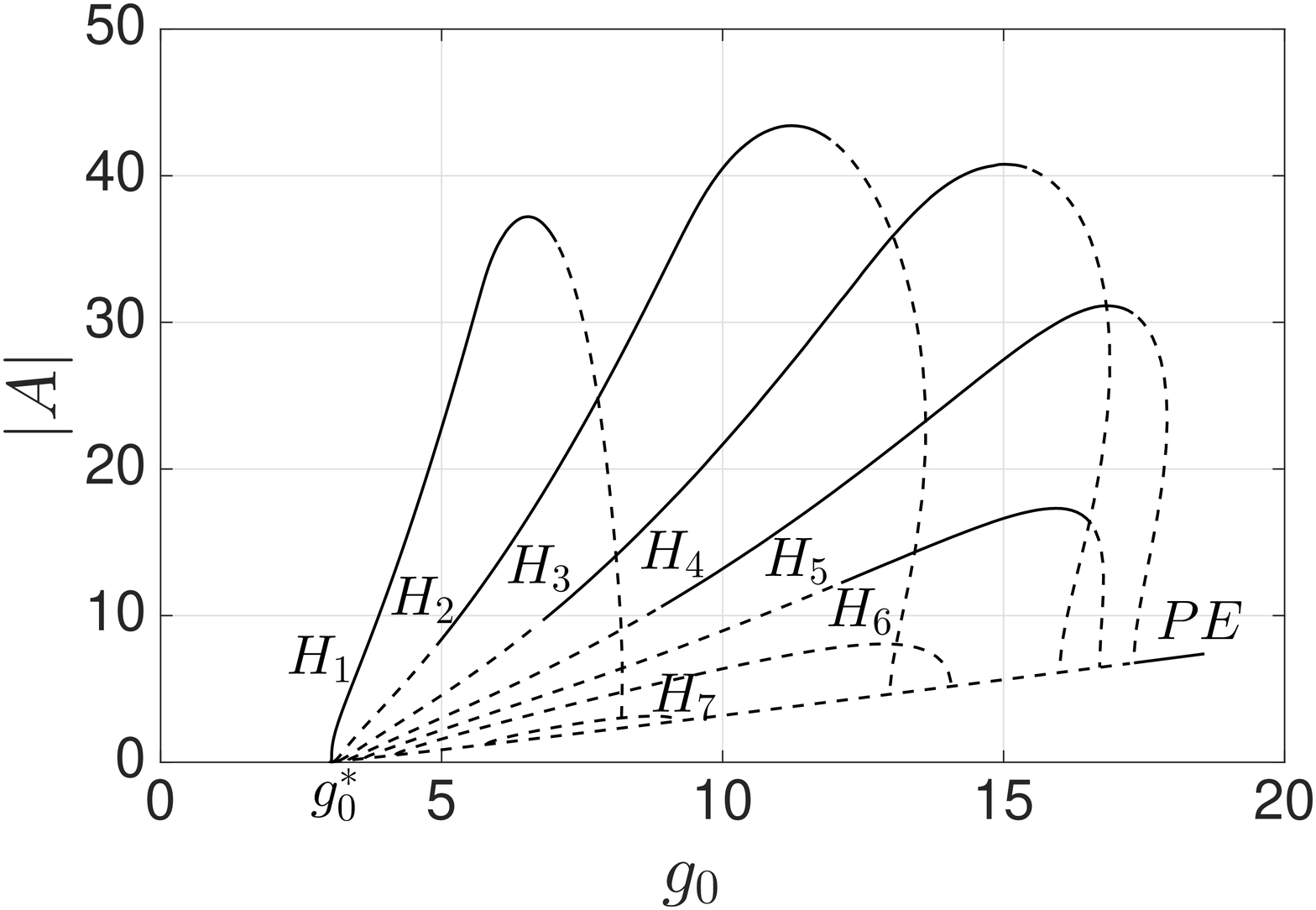}}
\hfill
\subfloat[\label{bifdiag3}]{\includegraphics*[width=.5\textwidth]{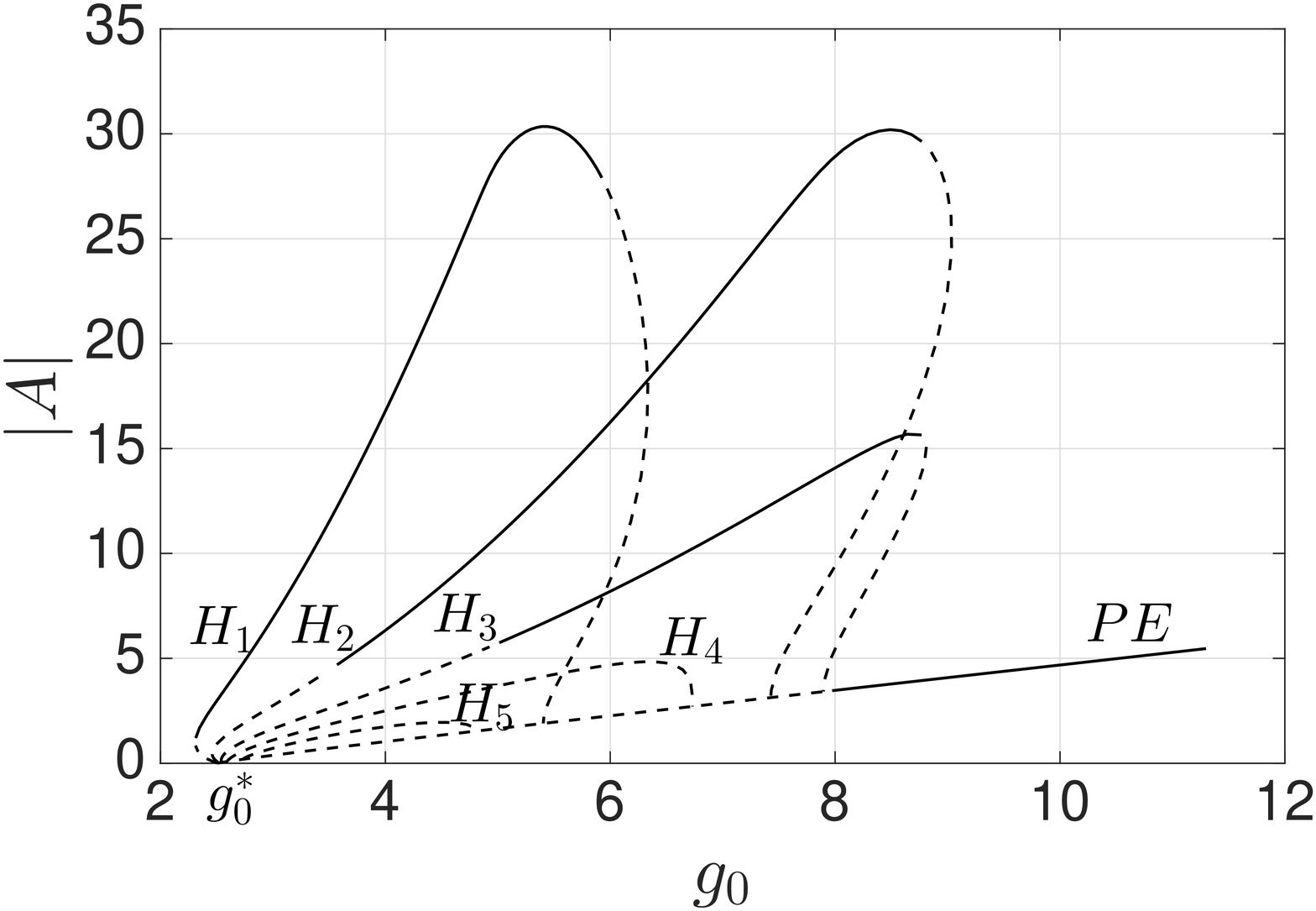}}
\end{center}
\caption{\small{Bifurcation diagrams obtained with numerical package DDE-BIFTOOL\cite{Engelborghs1, sieber2014dde, engelborghs2001dde, janssens2010normalization, wage2014normal} for system \eqref{A}--\eqref{G} for two parameter sets. The vertical axis shows the maximum of the $A$-component of a periodic solution. The $PE$ line corresponds to the positive equilibrium. Branches $H_1$--$H_7$ ($H_1$--$H_5$ on panel (b)) correspond to the periodic solutions born via Hopf bifurcations on the positive equilibrium. Stable branches are shown by solid lines and unstable branches are shown by dashed lines. The branch $H_1$ on panel (b) exhibits slight hysteresis near the threshold $g_0^*$. All the branches connect to the branch of the positive equilibrium at Hopf bifurcation points at both ends.}}
\label{figDiag}
\end{figure}

\begin{figure}[ht]
\begin{center}
%
\subfloat[\label{fig1a}]{\includegraphics*[width=.5\textwidth]{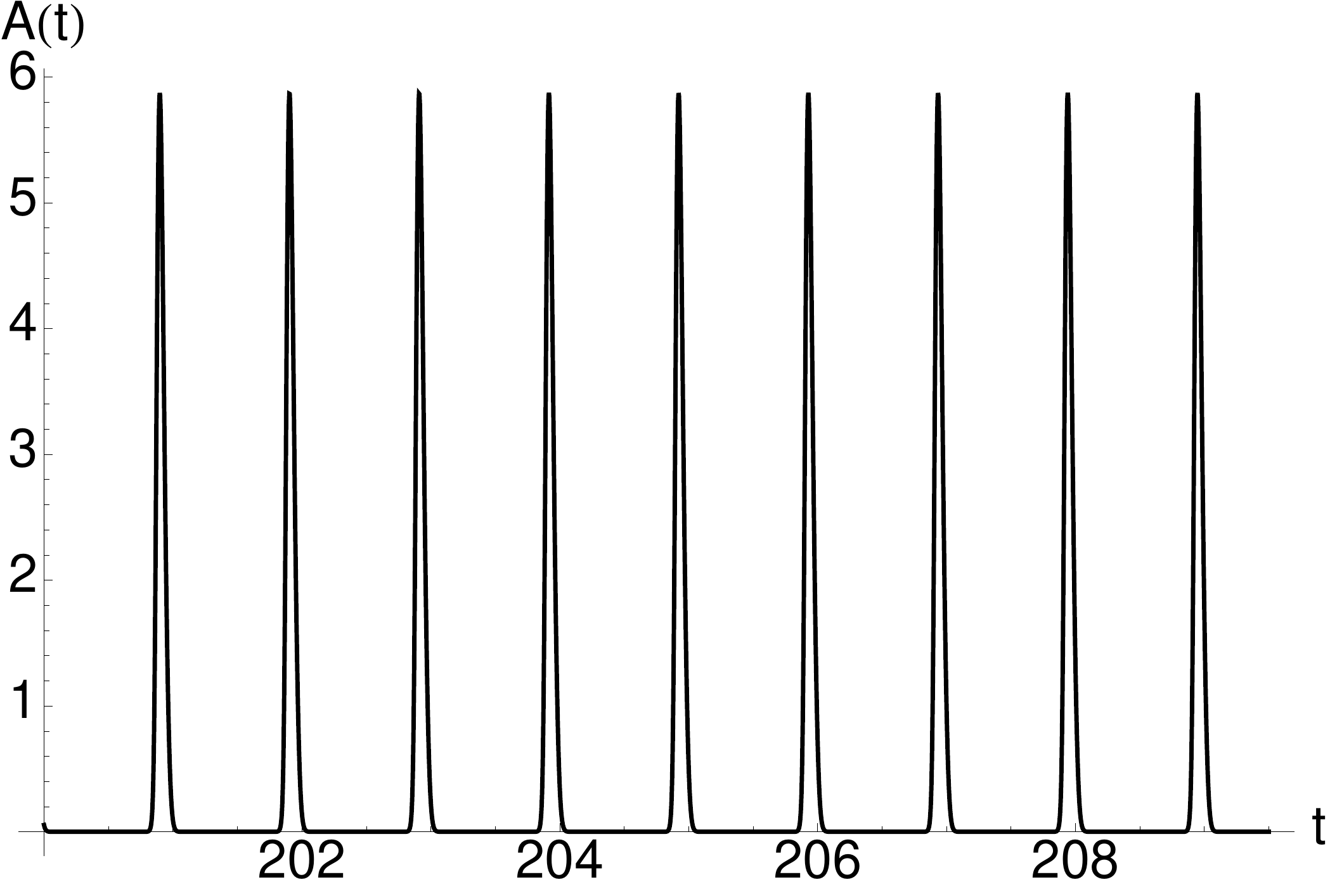}}
\hfill
\subfloat[\label{fig1q}]{\includegraphics*[width=.5\textwidth]{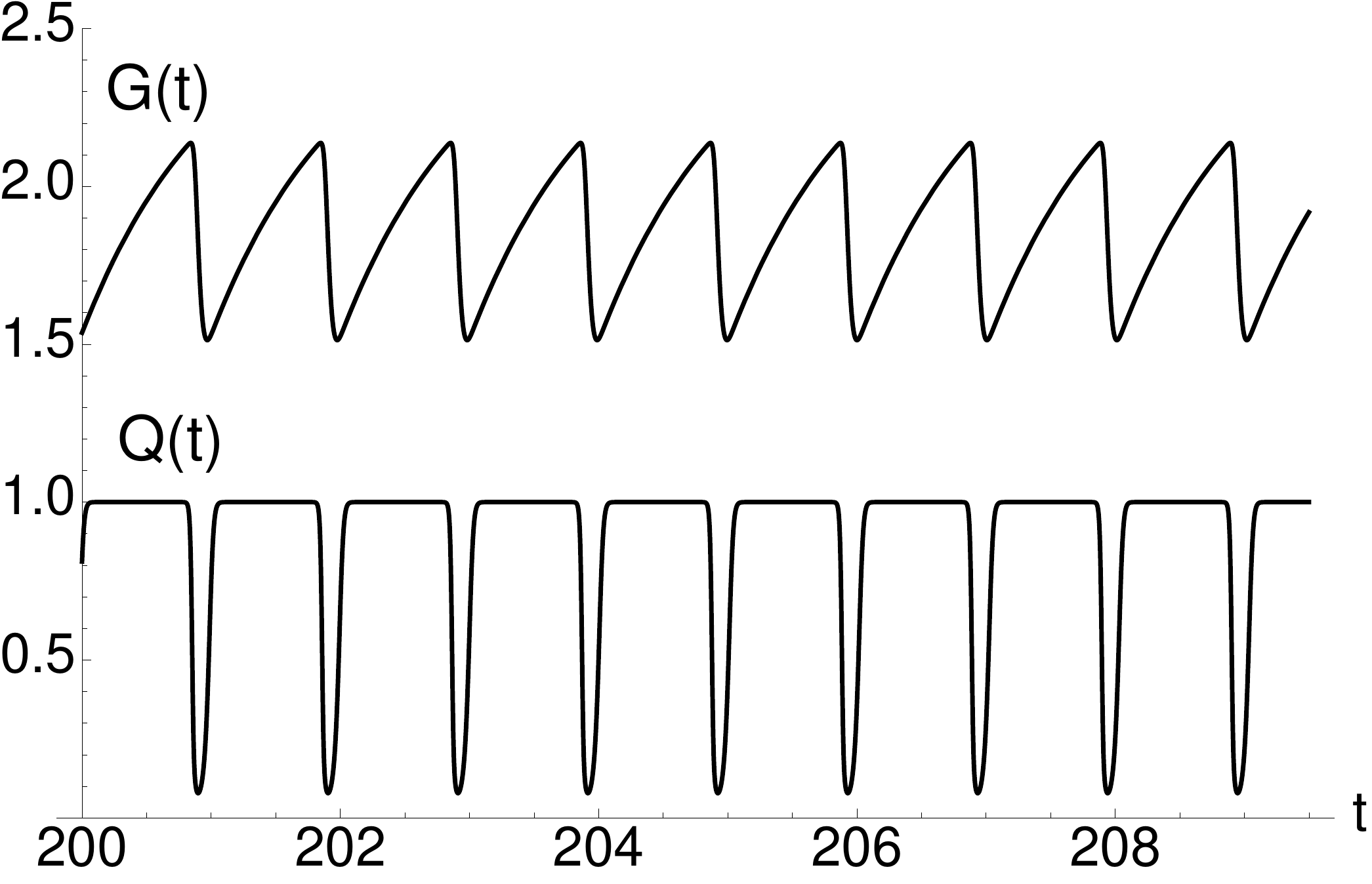}}
\end{center}
\caption{\small{Time trace of the periodic solution of system \eqref{A}--\eqref{G}. Panel (a): the $A$-component; Panel (b): the $G$-component (above) and the $Q$-component below.
The $A$-component is almost zero between the pulses.
The $Q$-component stays close to the equilibrium value $q_0/\beta=1$ between the pulses of the $A$-component and drops almost to zero during the pulse.
The $G$-component drops fast during the pulse and then recovers slowly between the pulses.
The period of the solution is close to the delay $T=1$. The following parameters were used: $\gamma= 200,\ \kappa=0.5,\ g_0= 3.1,\ q_0= 1,\ \alpha= 1,\ \beta= 1,\
s = 1,\ k = 1,\ T = 1$. The threshold value is $g_0^*=3$. \label{main_trace}}}
\end{figure}

\begin{figure}[ht]
\begin{center}
%
\subfloat[\label{phases}]{\includegraphics*[width=.5\textwidth]{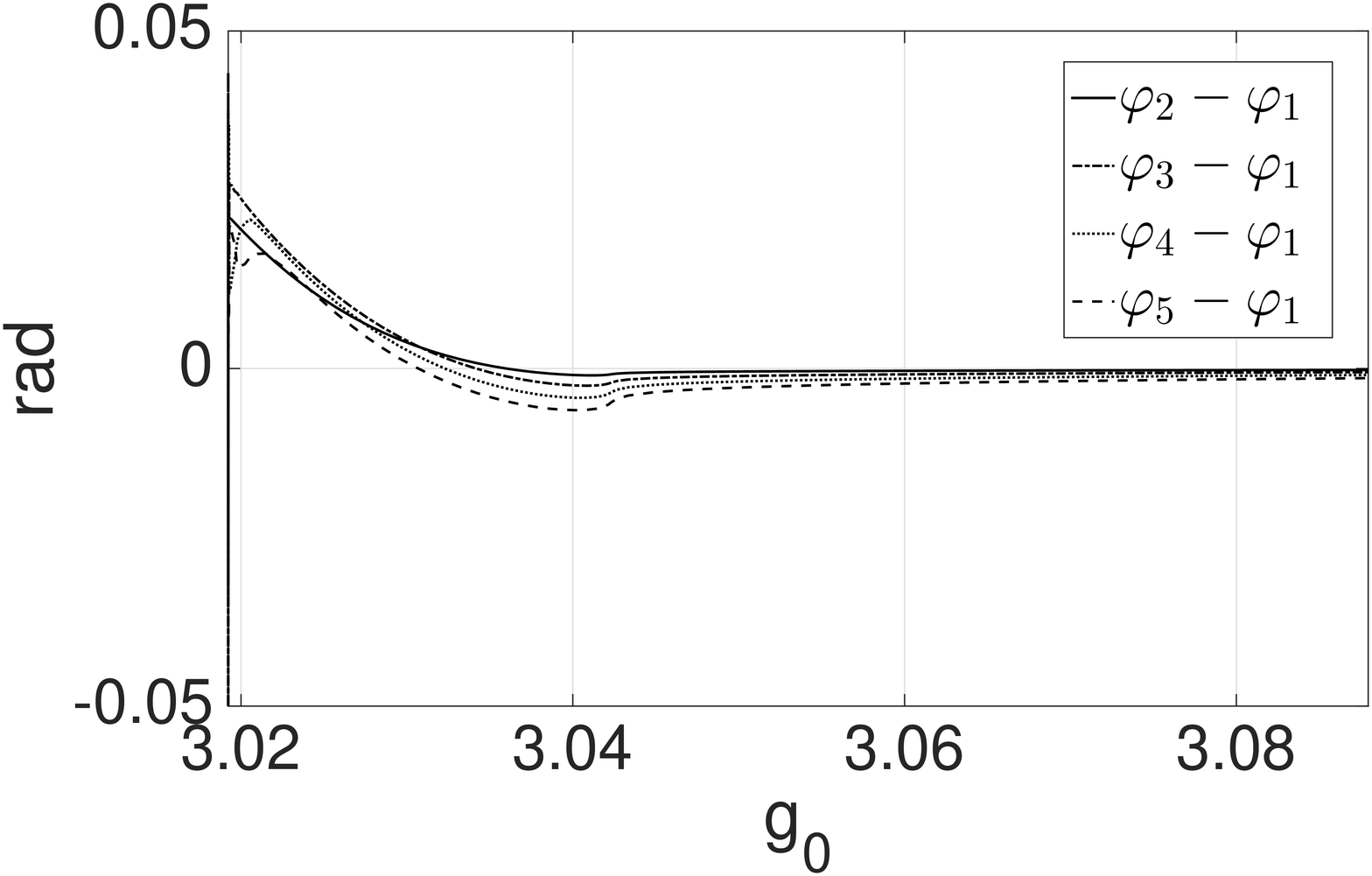}}
\hfill
\subfloat[\label{ampls}]{\includegraphics*[width=.5\textwidth]{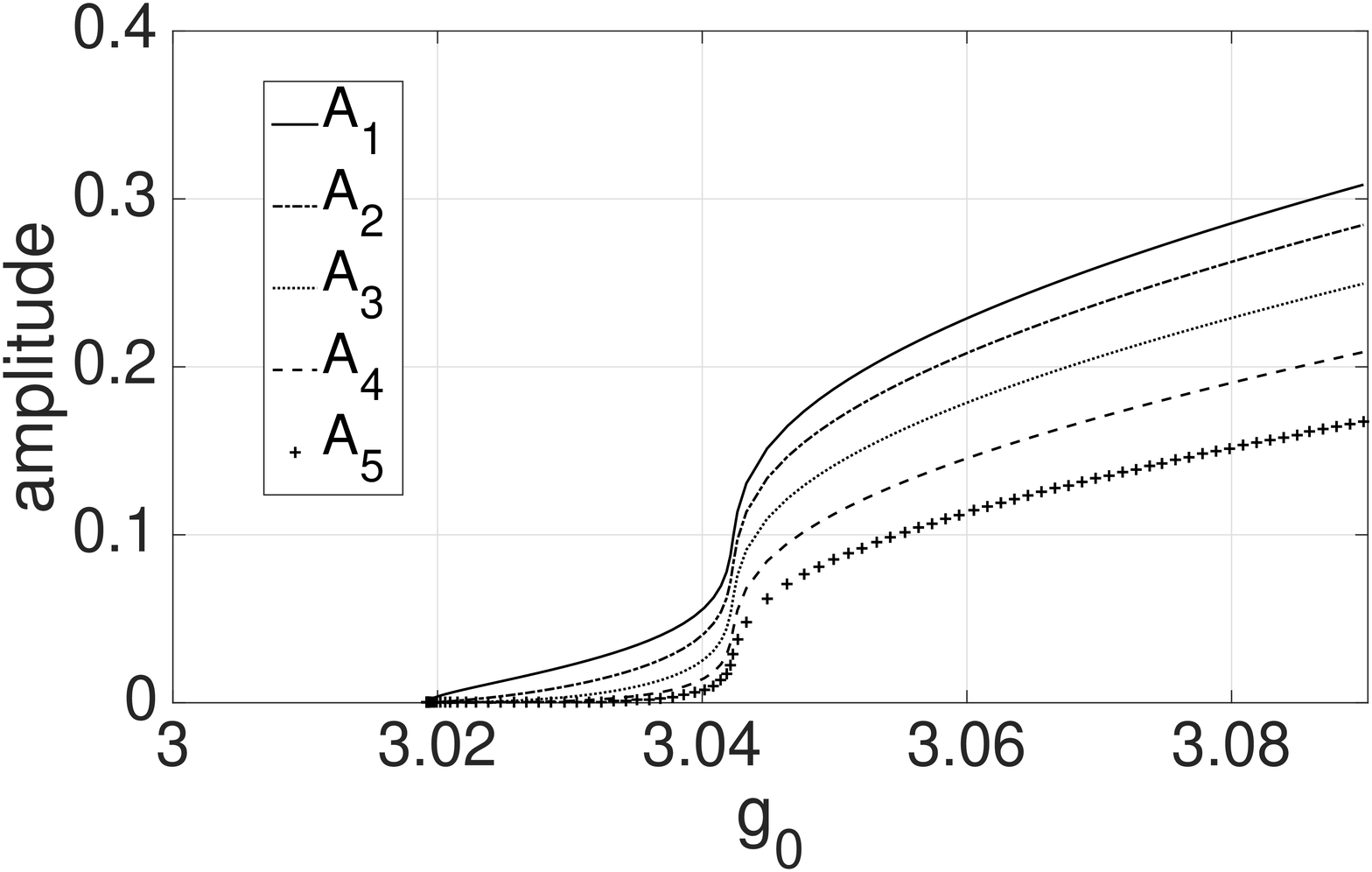}}
\end{center}
\caption{The phase and amplitude of the Fourier coefficients for the $A$-component $A(t/\tau) = \sum_{n=1}^\infty {A_n \cos{(2 \pi n t/\tau + \phi_n)}}$ of the periodic solution along the branch $H_1$ shown in Figure \ref{bifdiag2} where $\tau$ is the period of solution. \label{fourier}}
\end{figure}

As the bifurcation parameter $g_0$ increases across the threshold, the positive equilibrium \eqref{equil2} also undergoes a sequence of Hopf bifurcations,
which we deem responsible for the creation and formation of the periodic pulsating solution.
The first Hopf bifurcation with the frequency close to $2\pi/T$ destabilizes the positive equilibrium and creates a stable
cycle {(see branch $H_1$ in Figure \ref{figDiag})}. As the parameter $g_0$ increases further, this cycle changes its shape continuously into a pulsating periodic solution, see Figure \ref{main_trace}.  The amplitudes of harmonics of the $A$-component $A(t/\tau) = \sum_{n=1}^\infty {A_n \cos{(2 \pi n t/\tau + \phi_n)}}$ of the periodic solution, where $\tau$ is the period of $A$, grow with $g_0$, while the phase differences $\phi_k-\phi_1$ almost vanish, see Figure \ref{fourier}. 
At the same time the positive equilibrium undergoes a cascade of the secondary Hopf bifurcations with the frequencies of the higher harmonics.
 The whole cascade of the Hopf bifurcations and the transformation of the cycle
 to a pulsating solution happen in a small right neighborhood of the threshold $g_0=g_0^*$, see Figure \ref{fourier}.

 The characteristic equation for the positive equilibrium is
 \begin{equation} \label{positiveChar}
 	e^{T\lambda} = \frac{G_* \kappa \left(\alpha + \lambda\right)\left(A_* s + \beta + \frac \lambda \gamma \right)} {\left(A_* k+\alpha + \lambda\right)\left(A_* s \left(1+\frac \lambda \gamma\right) + \left(\beta + \frac \lambda \gamma\right)\left(1+ \frac \lambda \gamma+Q_* \mu\right)\right)}.
 \end{equation}
 Using asymptotic formulas \eqref{equil2} 
 and the ansatz $\lambda=i \omega$ for the eigenvalues of the linearization, we obtain the following asymptotic formulas for
 the frequency and the bifurcation value of the parameter at each Hopf bifurcation point:
 \begin{equation}\label{omega1}
 \fl \omega_n = \frac{2 \pi n}{T}\left(1-\frac{\alpha}{\kappa g_0^* \gamma T} + \frac{\alpha^2}{\left(\kappa g_0^*\gamma T \right)^2} \left(\frac{\frac{\beta ^2 \kappa g^*_0 k \left(2 \pi ^2 n^2-\alpha  T\right)}{\alpha T \left(\alpha^2+\left(\frac{2 \pi  n}{T}\right)^2\right)}+\mu  q_0 s}{\mu  q_0 s -\frac{\beta ^2 \kappa g^*_0 k}{\alpha^2+\left(\frac{2 \pi  n}{T}\right)^2}}\right)\right)+ O \left(\gamma^{-3} \right),
 \end{equation}
 \begin{equation}\label{delta1}
 \delta_n =\left(\frac{2 \pi n}{\gamma T}\right)^2\frac{\beta ^2 g^*_0 \kappa  k - \alpha ^2 \mu  q_0
    s}{2 \kappa^2 g^*_0 \left(\mu  q_0 s-\frac{\beta ^2
    g^*_0 \kappa  k} {\alpha ^2+\left(\frac{2 \pi n}{T}\right)^2 } \right)} + O \left(\gamma^{-3} \right).
 \end{equation}
  We assume that, along with the relation \eqref{trans}, the condition
  \begin{equation}\label{trans'}
\mu  q_0 s>\frac{\beta ^2
    g^*_0 \kappa  k} {\alpha ^2+\left(\frac{2 \pi}{T}\right)^2 }
\end{equation}
is satisfied. Under this condition, relation \eqref{omega1} implies
$\delta=g_0-g_0^*>0$ for $n=1,2,\ldots$\ That is, according to Eqs. \eqref{trans}, \eqref{omega1},
condition \eqref{trans'} ensures that the positive equilibrium undergoes the Hopf
bifurcations with the frequencies close to the multiples $2\pi n/T$ of $2\pi /T$ for
$n=1,2,\ldots$ as $g_0$ increases across the threshold.

\begin{figure}[ht]
\begin{center}
\subfloat[\label{stabil11}]{\includegraphics*[width=.5\textwidth]{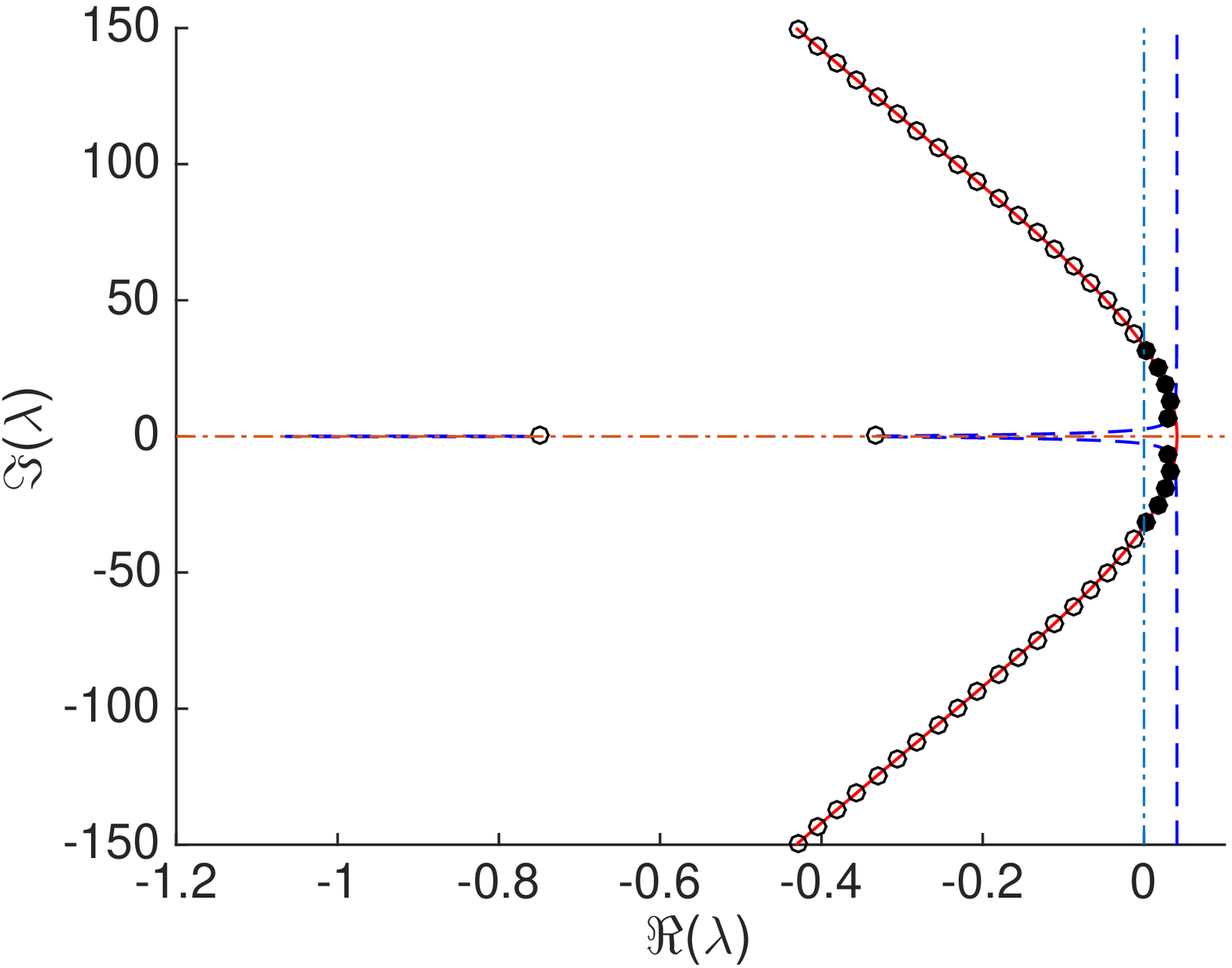}}
\hfill
\subfloat[\label{stabil5}]{\includegraphics*[width=.5\textwidth]{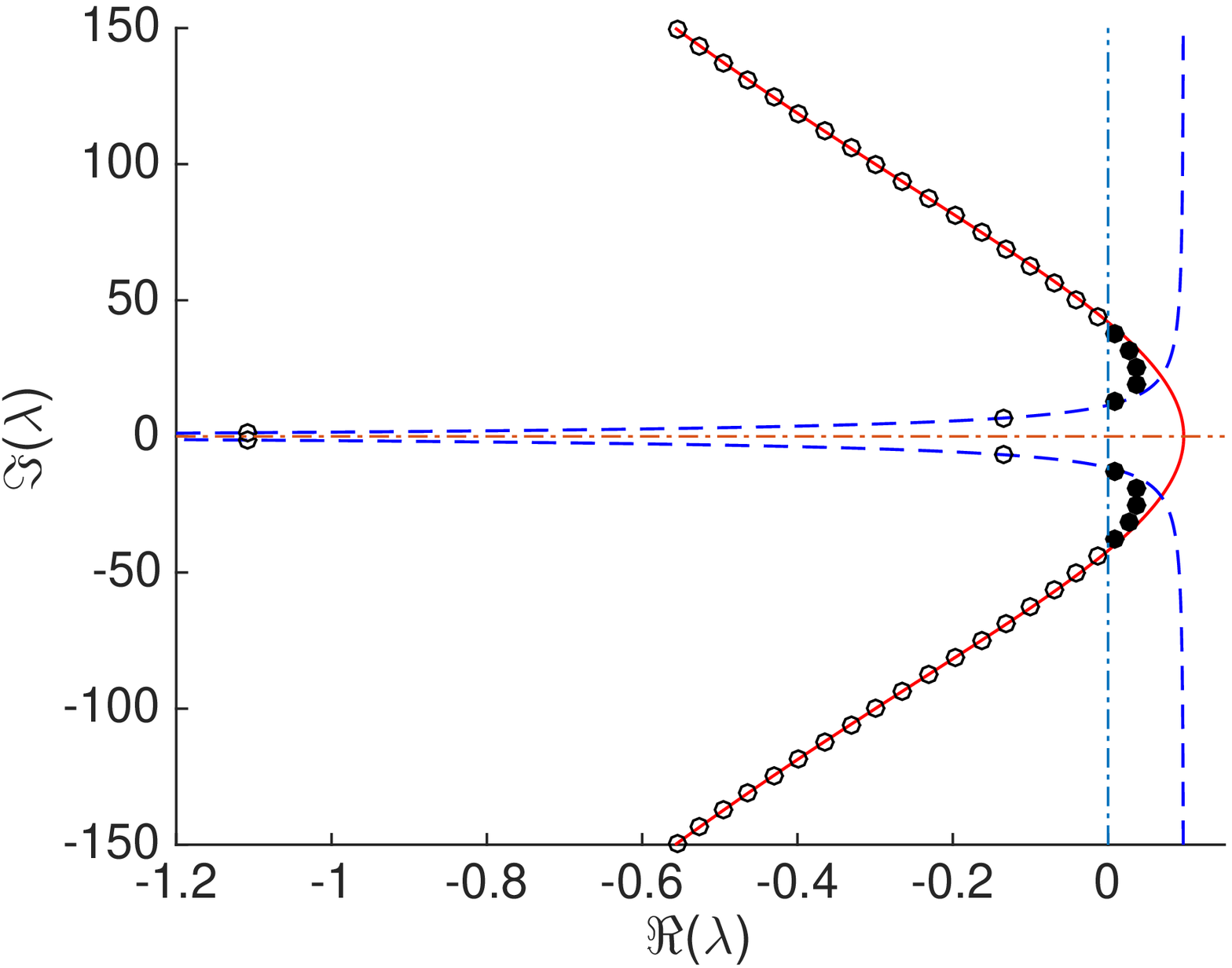}}
\end{center}
\caption{\small{Spectrum of the positive equilibrium and curves \eqref{weak}, \eqref{strong} for $g_0 = 3.7437$ (panel (a)) and $g_0 = 12.5436$ (panel(b)).
Notation and other parameters are the same as in Figure \ref{stabil}b.}}
\label{stabil1}
\end{figure}

The spectrum of the positive equilibrium can be divided into two parts, which have different asymptotic properties
with respect to the large parameter $\gamma$, cf.~\cite{wy}.
{\sl Strong} spectrum consists of the eigenvalues $\lambda = x+i\omega+O(\gamma^{-1})$, which originate from the limit $\gamma = \infty$.
\eqref{positiveChar} implies the following approximate implicit relationship
between the real and imaginary parts for these eigenvalues:
\begin{equation}\label{strong}
\frac{G_*^2 \kappa ^2 \left(A_* s+\beta \right){}^2 \left((\alpha +x)^2+\omega
   ^2\right)}{\left( (A_* k+\alpha +x)^2+\omega
   ^2\right) \left(A_* s+\beta +\beta  \mu  Q_*\right){}^2}-e^{2 T x} =0.
\end{equation}
\textit{Weak} spectrum is characterized by the asymptotic relationship $\lambda = x+i \gamma \omega$ and satisfies the approximate relationship
\begin{equation}\label{weak}
x(\omega) = \frac{1}{2T}\ln \left(\frac{G_*^2 \kappa ^2 \left(\left(A_*
   s+\beta \right){}^2+\omega ^2\right)}{ \left(A_* s+\beta +\beta  \mu
   Q_*\right){}^2+ W\omega ^2 +\omega ^4}\right),
\end{equation}
where
\begin{equation*}
	W = 2 A_* s
	   \left(\beta +\mu  Q_*\right)+A_*^2 s^2+\beta^2+\left(\mu
	   Q_*+1\right)^2.
\end{equation*}
With increasing $g_0$, the curve \eqref{weak} that carries the weak spectrum moves to the right
producing the Hopf bifurcations described by Eqs.~\eqref{omega1}, \eqref{delta1},
see Figures \ref{stabil}b and \ref{stabil1}a.
However, for larger values of $g_0$, weak eigenvalues with smaller imaginary part
leave this curve, become a part of the strong spectrum, and stabilize, see Figure \ref{stabil1}b.


In Tables \ref{tab:delta90}, \ref{tab:delta400}, the asymptotic values of $\omega_n$ and $\delta_n$ given by formulas \eqref{omega1} and \eqref{delta1} are compared with the numerical values obtained for two values of the parameter $\gamma$. Table \ref{tab:delta90} corresponds to branches $H_1$--$H_5$ in Figure \ref{bifdiag2}. Table \ref{tab:delta400} is obtained for a larger value of $\gamma=400$ and the same values of the other parameters. For $\gamma=90$ (Table \ref{tab:delta90}) we observe 7 branches of periodic solutions, for $\gamma=400$ (Table \ref{tab:delta400}) the number of branches increases to 24. The accuracy of the asymptotic formulas increases with decreasing $n$ and with increasing $\gamma$.

We have conducted a number of further numerical simulations with different parameter sets satisfying conditions \eqref{trans} and \eqref{trans'} and observed bifurcation diagrams and oscillating periodic solutions similar to those presented in Figures \ref{figDiag}, \ref{main_trace}.

\begin{table}[H]
\centering
\begin{tabular}{|c|c|c|c|c|c|c|}
\hline
 \multirow{2}{*}{\#}  & \multicolumn{3}{|c|}{$\delta_n=g_0-g_0^*$} & \multicolumn{3}{|c|}{$\omega_n$} \\ \cline{2-7}
 & Asymptotic & Numerical & Error (\%) & Asymptotic &  Numerical & Error (\%)\\ \hline
$H_1$  & 0.0191   & 0.0192  & 0.8 & 6.2382  & 6.2380 & 0.0026\\ \hline
$H_2$  & 0.0675   & 0.0712  &   5.1  & 12.4761 &  12.4750 & 0.0087 \\ \hline
$H_3$  & 0.1487   & 0.1701 &  12.5 & 18.7141 &  18.7105 & 0.0195 \\ \hline
$H_4$  & 0.2624   & 0.3391  &   22.6 & 24.9521 & 24.9432 & 0.0359 \\ \hline
$H_5$  & 0.4086   & 0.6296  &   35.1 & 31.1901 & 31.1716 & 0.0596\\ \hline
\end{tabular}
\caption{Comparison of the asymptotic and numerical values of $\delta_n$, $\omega_n$ for the parameter set used in Figure \ref{bifdiag2}.}
\label{tab:delta90}
\end{table}

\begin{table}[H]
\centering
\begin{tabular}{|c|c|c|c|c|c|c|}
\hline
 \multirow{2}{*}{\#}  & \multicolumn{3}{|c|}{$\delta_n = g_0 - g_0^*$} & \multicolumn{3}{|c|}{$\omega_n$} \\ \cline{2-7}
   & Asymptotic & Numerical &  Error (\%) & Asymptotic & Numerical &  Error (\%) \\ \hline
$H_1$  & $9.6559 \cdot 10^{-4}$   & $9.64 \cdot10^{-4}$   &   0.17 & 6.2728 & 6.2728 &  $4.8316\cdot10^{-6}$ \\ \hline
$H_2$  & 0.0034\ignore{0.003418955407243}  & 0.0034\ignore{0.003428496859645}  &  0.28  &  12.5456 & 12.5456 &  $8.0594\cdot10^{-5}$ \\ \hline
$H_3$  &  0.0075  & 0.0076  &  0.52   & 18.8184 & 18.8183 &  $1.3973\cdot10^{-4}$  \\ \hline
$H_4$  &  0.0133  &   0.0134   &   0.95  &  25.0911 & 25.0910 &  $4.1168\cdot10^{-4}$  \\ \hline
$H_5$  &  0.0207  &   0.0210   &    1.65   &  31.3639 & 31.3637 &  $6.4092\cdot10^{-4}$ \\ \hline
$H_{10}$  & 0.0824   & 0.0890   &    7.44   &  62.7278 & 62.7262 &  0.0026 \\ \hline
$H_{20}$  &  0.3291   &  0.4654  &   29.29  & 125.4557 & 125.4409 &  0.0118 \\ \hline
\end{tabular}
\caption{Comparison of the asymptotic and numerical values of $\delta_n$, $\omega_n$ for $\gamma=400$ and other parameters are the same as in Table \ref{tab:delta90}.}
\label{tab:delta400}
\end{table}

\FloatBarrier


\section{Variations of the main model}\label{variations}
We now consider three variations of the model Eqs. \eqref{A}--\eqref{G}.
First, we consider the system of equations \eqref{A}, \eqref{G} and
\begin{equation}\label{Q'}
\gamma^{-1}Q'=q_0 Q-\beta Q^2-s AQ.
\end{equation}
That is, we replace the constant immigration term and the linear death term in the $Q$-equation
by the logistic growth term $Q(q_0-\beta Q)$.
System \eqref{A}, \eqref{G}, \eqref{Q'} has the same
equilibrium \eqref{equil1} as system \eqref{A}--\eqref{Q}, and this equilibrium undergoes
the transcritical bifurcation at the same threshold $g_0=g_0^*$ defined by relation
\eqref{thresho}. Assuming the relation
\begin{equation}\label{trans1}
\frac{ k g_0^*}{\alpha^2}>\frac{s \mu}{\beta \kappa},
\end{equation}
which is a counterpart of {relation} \eqref{trans}, we ensure that the positive equilibrium
exists for $g_0>g_0^*$ and is stable in a right neighborhood of the threshold.
The components of the positive equilibrium near the threshold are defined by
the asymptotic formulas \eqref{equil2} where the coefficients now have the form
\begin{equation}
\tilde a =\frac 1 {\frac{k g_0^*}\alpha-\frac{\alpha s \mu}{\beta \kappa}},\qquad \tilde q=
\frac 1{ \frac {\mu \alpha} {\kappa}-\frac{k g_0^*\beta}{\alpha s}}, \qquad \tilde g = \frac \mu \kappa \tilde q.
\end{equation}
System \eqref{A}, \eqref{G}, \eqref{Q'}
has an additional equilibrium point $A=Q=0$, $G=g_0/\alpha$, which is unstable
for all $g_0>0$.

\begin{figure}[ht!]
\begin{center}
\subfloat[]{\includegraphics*[width=.5\columnwidth]{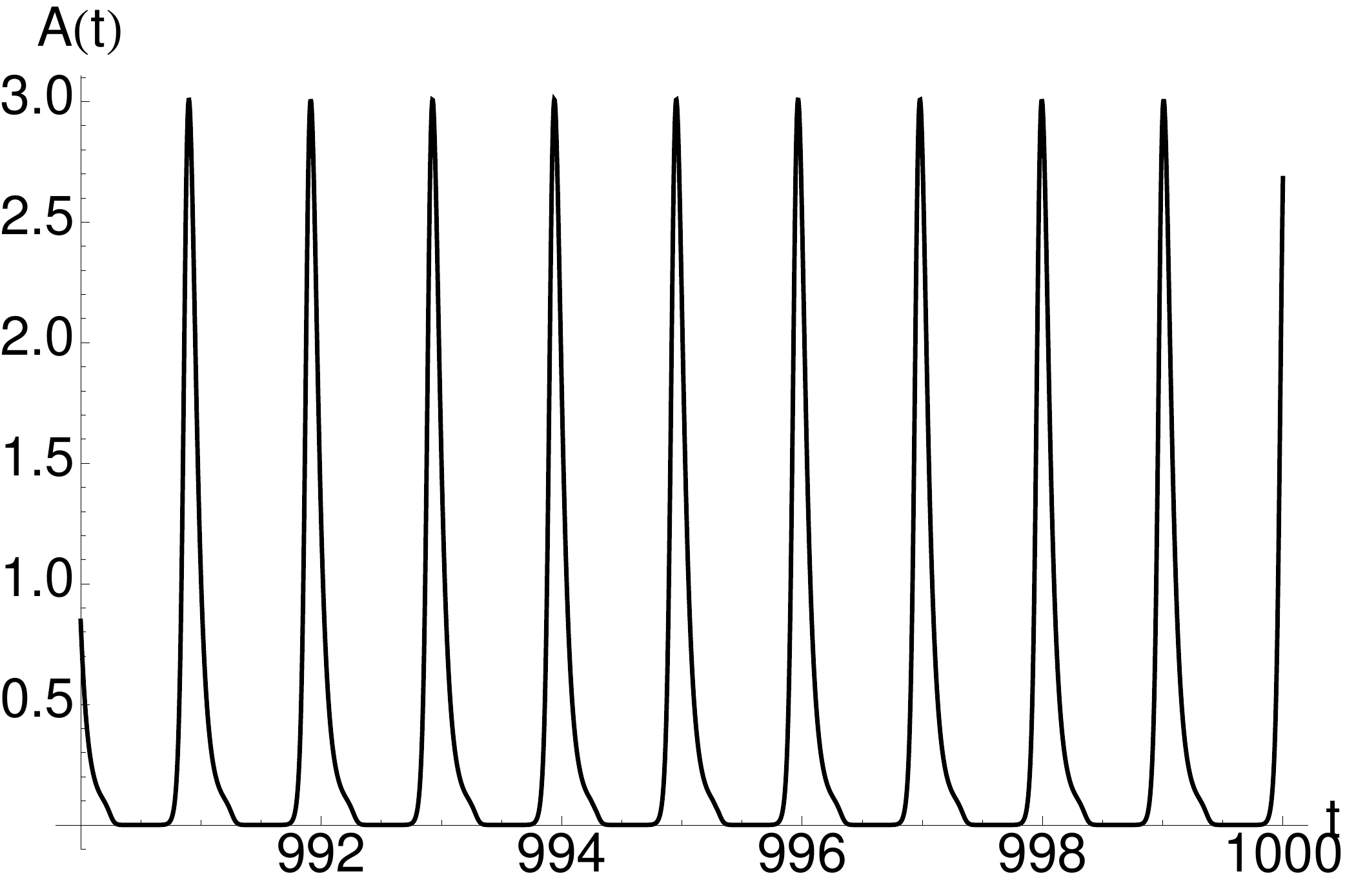}}
\hfill
\subfloat[]{\includegraphics*[width=.5\columnwidth]{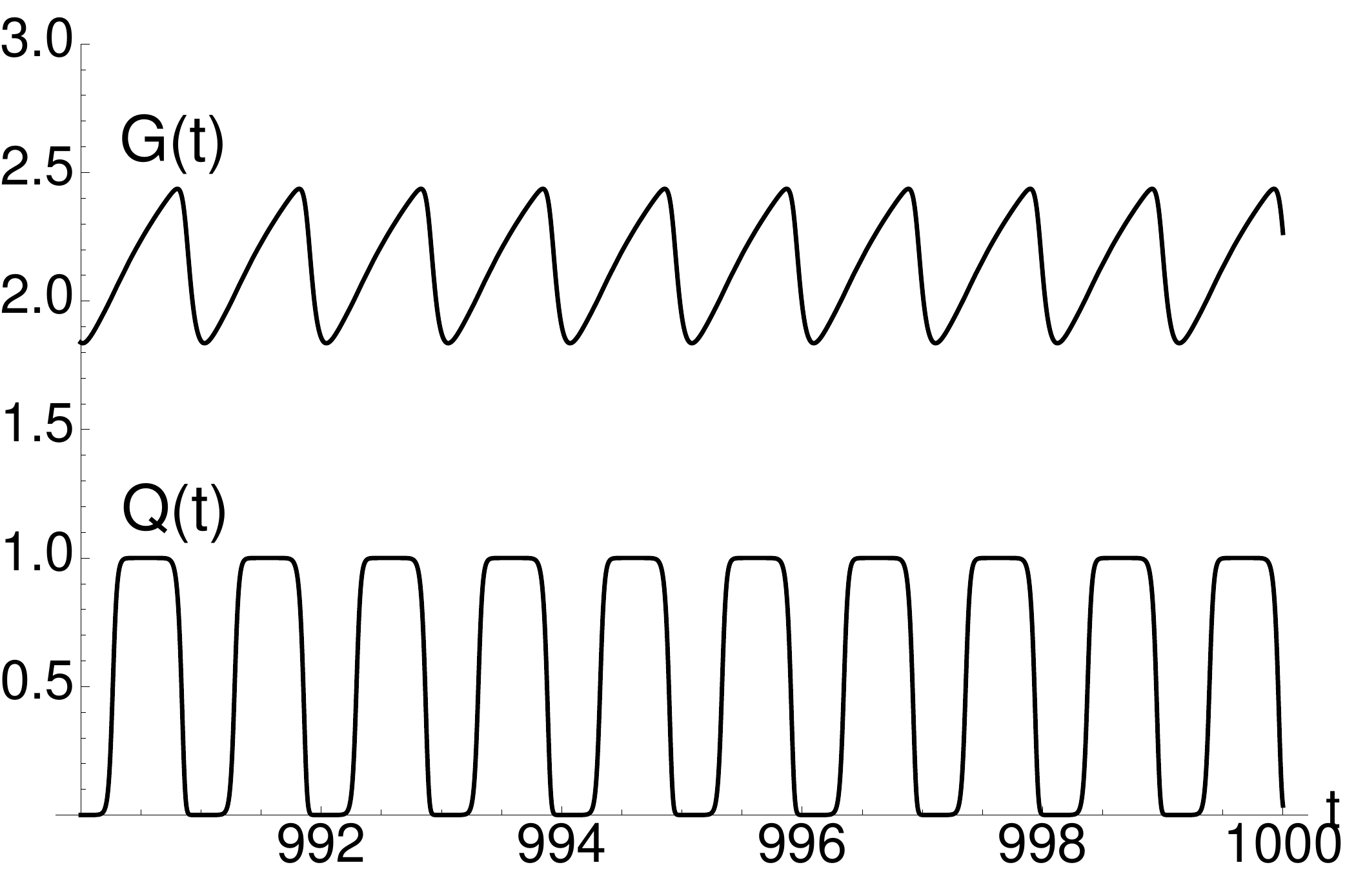}}
\end{center}
\caption{Time traces of the $A$-component (panel (a)), the $G$-component and the $Q$-component (panel (b))
of the cycle of Eqs. \eqref{A}, \eqref{G}, \eqref{Q'}. The picture is similar to that of the cycle of system Eqs. \eqref{A}--\eqref{G}
in Figure \ref{main_trace}. Here $\gamma=60$ and $g_0=3.05$; the other parameters and the threshold are the same as in Figure \ref{main_trace}. \label{fig4}}
\end{figure}

The linearization {at} the equilibrium \eqref{equil1} with $A=0$ leads to the same
characteristic equation \eqref{eig} for the eigenvalues. Hence, this equilibrium looses stability at the threshold
and then undergoes a cascade of Hopf bifurcations at the same points \eqref{bifi:delta} and with the same frequencies \eqref{bifi}
as model \eqref{A}--\eqref{G} with increasing the bifurcation parameter $g_0$.
If, instead of \eqref{trans'}, we assume the relation
\begin{equation}\label{trans2}
\frac{ k g_0^*}{\alpha^2+\left(\frac{2\pi}T\right)^2}< \frac{s \mu}{\beta \kappa},
\end{equation}
which is simultaneous with \eqref{trans1}, then the positive equilibrium
undergoes a cascade of Hopf bifurcations in a small right neighborhood of the threshold.
The frequencies and the bifurcation values of the parameter $\delta=g_0-g_0^*$ for these bifurcations are
approximated by the formulas
\begin{equation*}
\fl \omega_n = \frac{2 \pi n }{T} \left(1-\frac{\alpha }{\kappa g_0^* \gamma T}+ \frac{\alpha^2 }{\left(\kappa g_0^* \gamma T\right)^2} \left( \frac{ \frac{ \beta \kappa g_0^* k
   \left(2 \pi^2 n^2-\alpha T\right)}{\alpha T \left(\alpha^2 + \left(\frac{2 \pi n}{T} \right)^2\right)}+\mu  s}{\mu s - \frac{k \beta \kappa g_0^*}{\alpha^2 + \left(\frac{2 \pi n}{T} \right)^2}}\right)\right)+O\left(\gamma ^{-3}\right),
\end{equation*}
\begin{equation*}
\delta_n = \left(\frac{2 \pi n}{\gamma T} \right)^2\frac{ \beta k g_0^* \kappa - \alpha^2  \mu  s}{2\kappa^2
   g_0^* \left(\mu  s - \frac{k \beta \kappa g_0^*}{\alpha^2 + \left(\frac{2 \pi n}{T} \right)^2}\right)}+O\left(\gamma^{-3}\right)
\end{equation*}
for $n=1,2,\ldots$\ Figure \ref{fig4} presents the stable cycle of system \eqref{A}, \eqref{G}, \eqref{Q'}.

\begin{figure}[ht!]
\begin{center}
\subfloat[]{\includegraphics*[width=.5\columnwidth]{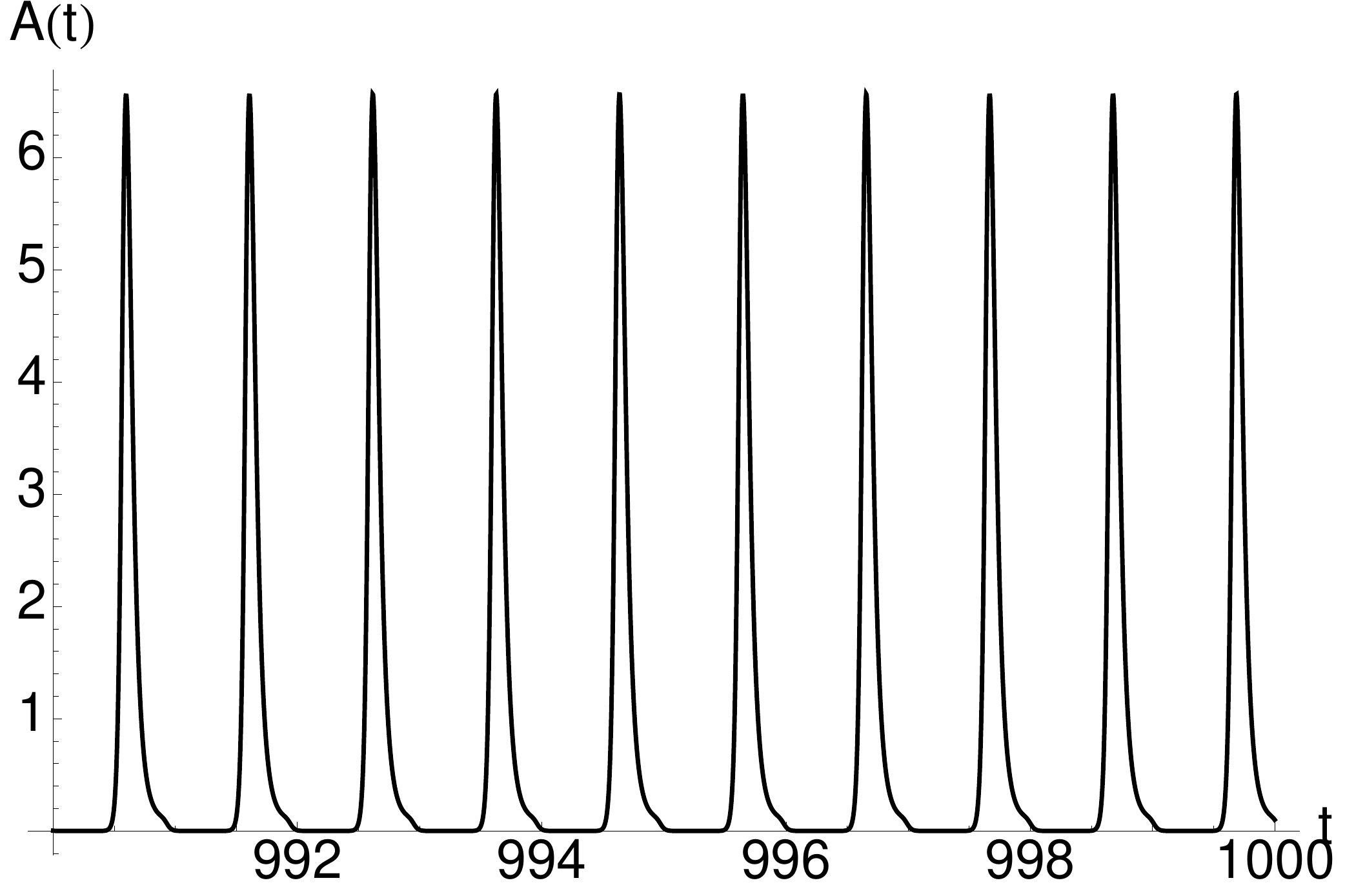}}
\hfill
\subfloat[]{\includegraphics*[width=.5\columnwidth]{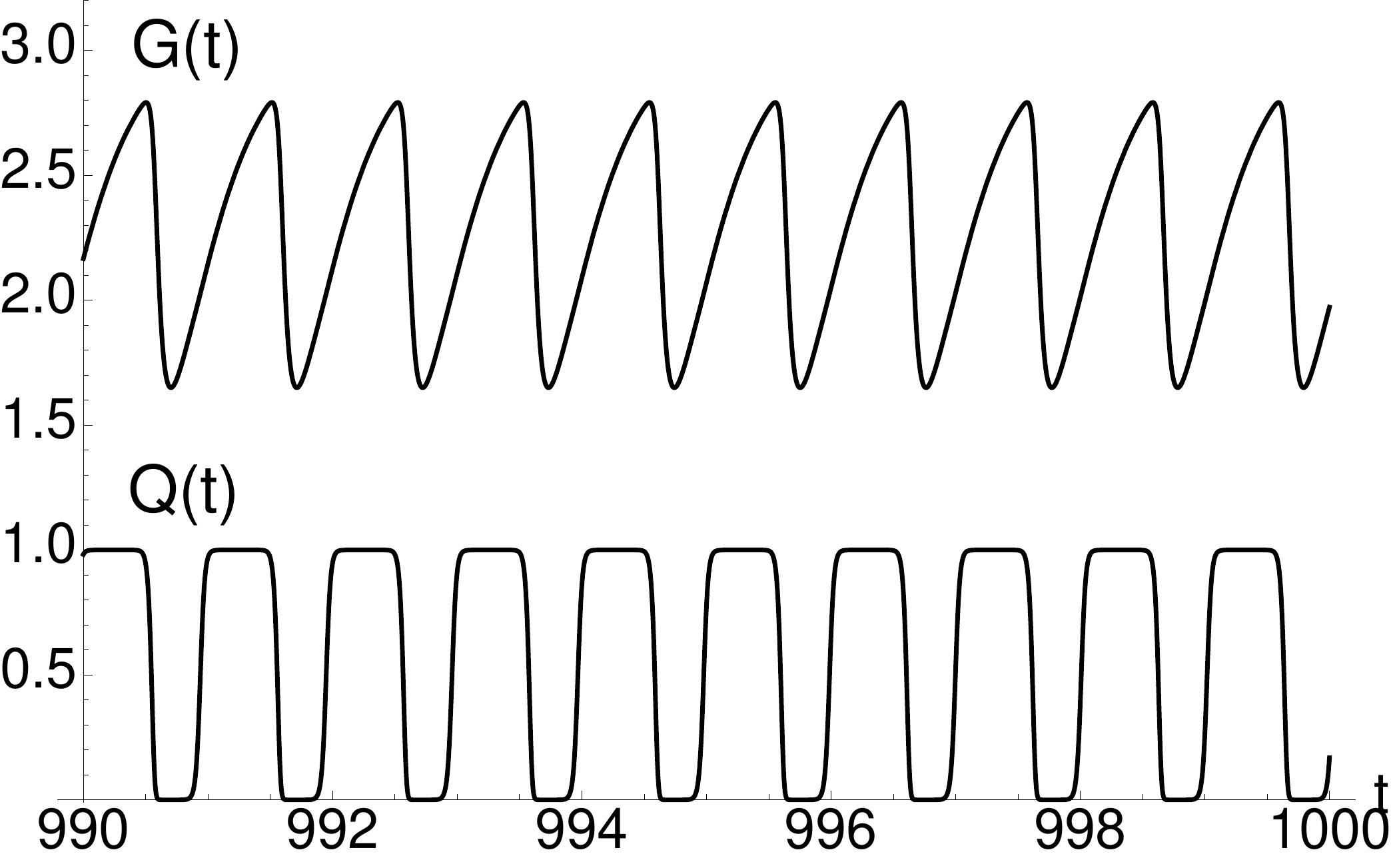}}
\end{center}
\caption{The stable cycle of Eqs. \eqref{AA}-\eqref{GG}. Here $s=0.5$; the other parameters and the threshold are the same as in Figure \ref{fig4}. \label{fig5}}
\end{figure}

As the next example, we consider the system of equations
\begin{eqnarray}
\gamma^{-1} A' &=& -A + \kappa G(t-T) A(t-T) - \mu Q A,\label{AA}\\
\gamma^{-1} Q' &=& q_0 Q-\beta Q^2 -s A Q, \label{QQ}\\
 G' &=& g_0 G-\alpha G^2 - k A G \label{GG}
\end{eqnarray}
where both the $Q$ and $G$ species have the logistic growth; the $A$-equation is the same as in the other examples.
As in the previous examples, the point \eqref{equil1} is an equilibrium of system \eqref{AA}--\eqref{GG},
which is stable for $g_0<g_0^*$ and unstable for $g_0>g_0^*$ with the threshold $g_0^*$
at the transcritical bifurcation point defined by relation \eqref{thresho}. The positive equilibrium,
which collides with equilibrium \eqref{equil1} at the threshold, is defined by the relations
\begin{equation}\label{eqpos}
\fl A=\frac{\beta \kappa \delta}{k\beta \kappa-s\alpha \mu},\qquad Q=\frac{q_0}\beta + \frac{\kappa s \delta}{\alpha \mu s-\kappa \beta k}, \qquad G=\frac{g_0^*}\alpha+\frac{\mu s \delta}{\alpha \mu s - \kappa \beta k}.
\end{equation}
We assume that
\begin{equation}\label{transs}
k\beta \kappa> s\alpha \mu,
\end{equation}
thus ensuring that the positive equilibrium exists for $g_0>g_0^*$ and is stable
in a right neighborhood of the threshold. The equilibrium \eqref{equil1} undergoes
the same cascade of Hopf bifurcations (with the same frequencies and at the same bifurcation points) above the threshold
as in the previous examples. The counterpart of condition \eqref{trans'} for Eqs. \eqref{AA}--\eqref{GG} is
\begin{equation}\label{trans''}
\alpha \mu s> \frac{ (g_0^*)^2 \kappa \beta k}{(g_0^*)^2+\left(\frac{2\pi}T\right)^2}.
\end{equation}
Relations \eqref{transs}, \eqref{trans''} imply that positive equilibrium \eqref{eqpos}
undergoes a cascade of Hopf bifurcations with the frequencies and bifurcation points
defined by
\begin{equation*}\label{omegaVar2}
\fl \omega = \frac{2 \pi n}{T}\left(1-\frac{\alpha}{\gamma T \kappa g_0^*}+\left(\frac{\alpha}{\gamma T \kappa g_0^*}\right)^2 \frac{\frac{\beta  g_0^* \kappa  k \left(2 \pi
   ^2 n^2-g_0^* T\right)}{T
   \left(\left(g_0^*\right)^2+\frac{4 \pi ^2 n^2}{T^2}\right)}+\alpha
   \mu  s}{\alpha  \mu
   s-\frac{(g_0^*)^2
   \kappa \beta k}{(g_0^*)^2+\left(\frac{2 \pi n}{T}\right)^2}}\right)+O(\gamma^{-3}),
\end{equation*}
\begin{equation*} \label{deltaVar2}
\delta = \frac{1}{2 g_0^*} \left(\frac{2 \pi n \alpha}{ \gamma T \kappa}\right)^2\frac{\kappa \beta k-\alpha  \mu  s}{\alpha  \mu  s-\frac{(g_0^*)^2
   \kappa \beta k}{(g_0^*)^2+\left(\frac{2 \pi n}{T}\right)^2}} + O(\gamma^{-3})
\end{equation*}
where $n=1,2,\ldots$\ Figure \ref{fig5} shows the stable cycle of system \eqref{AA}-\eqref{GG}.

{Finally, we consider the following model with competing fast species:}
\begin{eqnarray}
\gamma^{-1} A' &=&\kappa G(t-T)A(t-T) -\tau A - \mu Q A - f A^2, \label{A1}\\
\gamma^{-1} Q' &=& \nu G Q - \beta Q - s A Q - r Q^2, \label{Q1}\\
 G' &=& g_0 - \alpha G - k A G - m Q G. \label{G1}
\end{eqnarray}
Unlike previous examples, here we use the natural death rate $\tau$ of species $A$ as a bifurcation parameter to show the existence of a pulsating periodic solution. System \eqref{A1}--\eqref{G1} has an equilibrium $A_*=0$, $Q_*$, $G_*$ defined by
\begin{equation}
\nu G_*- r Q_*=\beta, \qquad \alpha G_* + m Q_* G_* = g_0.
\label{zeq}
\end{equation}
We assume
\begin{equation}
\nu g_0 > \alpha \beta.
\end{equation} This relation ensures that $G_*, Q_*>0$. System \eqref{A1}--\eqref{G1} also has a positive equilibrium which collides with the equilibrium \eqref{zeq} in a transcritical bifurcation for threshold value $\tau=\tau^*$ defined by
\begin{equation}\label{thresholdTau}
\tau^*= \kappa G_* - \mu Q_*.
\end{equation}

\begin{figure}[ht]
\begin{center}
\subfloat[]{\includegraphics*[width=.5\textwidth]{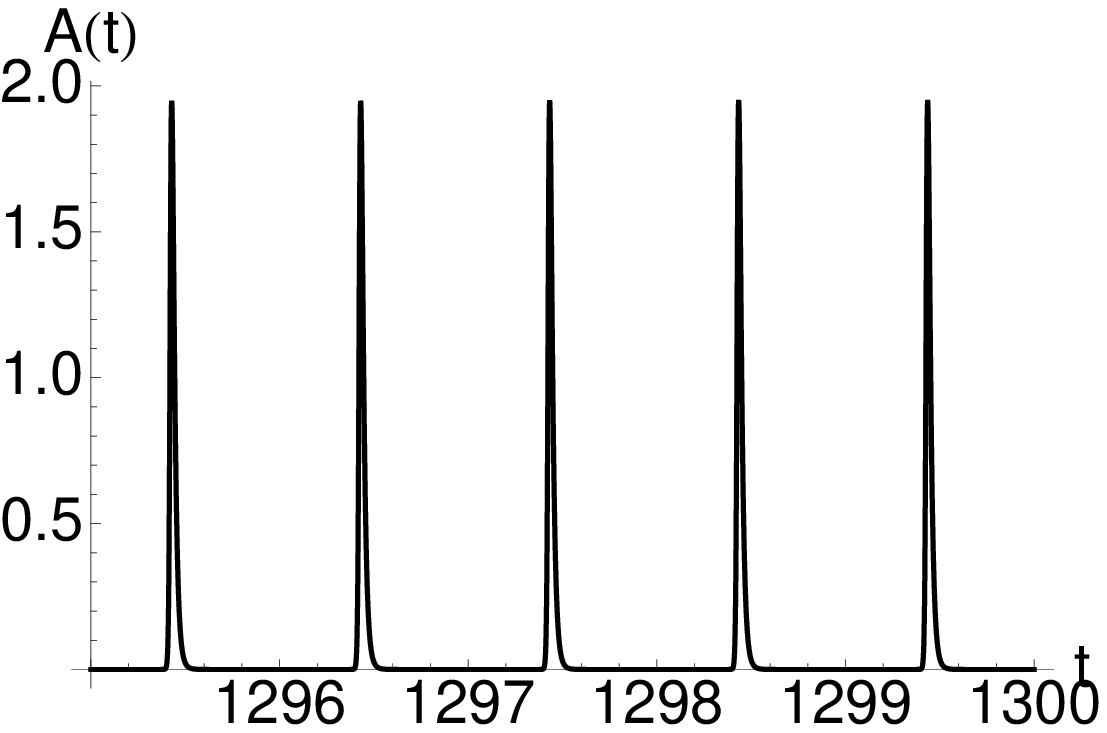}}
\hfill
\subfloat[]{\includegraphics*[width=.5\textwidth]{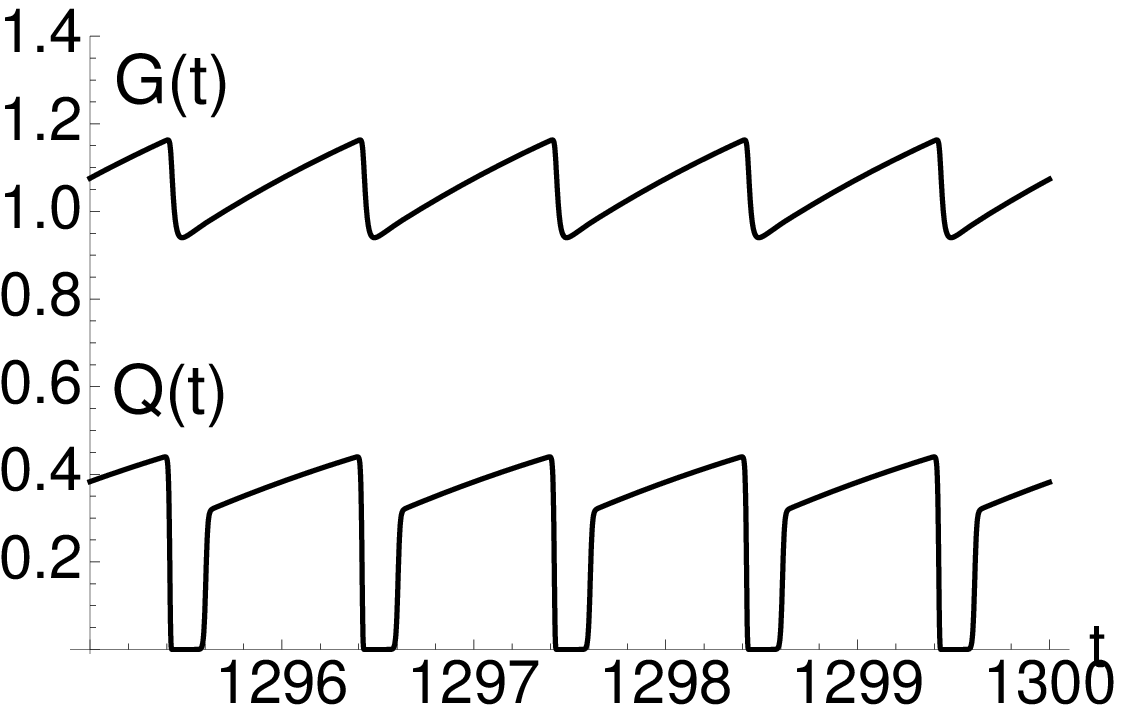}}
\end{center}
\caption{Time trace of solutions for the system \eqref{A1}--\eqref{G1}. The following parameters were used: $\gamma = 200, T= 1, \nu = 2, r= 3, \alpha = 0.3, m= 0.1,\kappa = 2, \mu = 1, s = 3, k = 4, f = 0.05, \beta = 1, g_0 = 0.6, \tau = 2$.}
\label{fig:fast:compete}
\end{figure}

The eigenvalues of the linearization of system \eqref{A1}--\eqref{G1} at the equilibrium \eqref{zeq} are defined by
\begin{eqnarray}
\label{zeroEigen1}\frac{\lambda ^2}{\gamma }+\lambda  \left(\frac{g_0}{\gamma  G_*}+Q_*
   r\right)+\frac{g_0 Q_* r}{G_*}+m \nu G_* Q_* &=&0, \\
\label{zeroEigen2}-\tau - \mu Q_* + \kappa G_* e^{-\lambda T} - \frac \lambda \gamma&=&0.
\end{eqnarray}
One can show that equations \eqref{zeroEigen1}--\eqref{zeroEigen2} imply $\mathrm{Re} \, \lambda < 0$ for $\tau>\tau^*$. For $\tau=\tau^*-\delta$ the components of the positive equilibrium near the threshold can be represented in the form of regular expansion in terms of $\delta$
\begin{equation}\label{poseq1}
\fl A = \tilde{A}\delta+O\left(\delta ^{2}\right), \qquad Q =Q_* + \tilde{Q}\delta+O\left(\delta ^{2}\right), \qquad G = G_* + \tilde{G}\delta+O\left(\delta ^{2}\right),
\end{equation}
where

\begin{equation*}
\tilde{A}=\frac{\frac{g_0 r}{G_*}+G_* m \nu }{F_*}, \qquad \tilde{Q} = -\frac{\frac{g_0 s}{G_*}+G_* k \nu }{F_*}, \qquad \tilde{G} = \frac{G_* (m s-k r)}{F_*},
\end{equation*}
with
\begin{equation*}
F_* = f \left(\frac{g_0 r}{G_*}+G_* m \nu
   \right)-\mu  \left(\frac{g_0 s}{G_*}+G_* k \nu \right)+\kappa G_* \left(k
   r- m s\right).
\end{equation*}
We assume that
\begin{equation}
F_* > 0
\end{equation} in order to guarantee that this equilibrium is positive for $\tau<\tau_*$ and hence the positive equilibrium is stable near the threshold. Furthermore we assume that
\begin{equation}\label{inter:intra}
\mu s > f r,
\end{equation}
and
\begin{equation}\label{last:cond}
F_* < \left(\frac{2 \pi}{T}\right)^2 \frac{r (\mu s-f r)}{\frac{g_0 r}{G_*}+G_* m \nu}.
\end{equation}
The conditions \eqref{inter:intra} and \eqref{last:cond} guarantee that the positive equilibrium undergoes the cascade of Hopf bifurcations in a small left neighborhood of the threshold $\tau^*$.
Condition \eqref{inter:intra} means that interspecific competition between fast species $A$ and $Q$ is stronger than intraspecific competition. In the classical competing species model this condition ensures the competive exclusion scenario; the opposite inequality $\mu s < f r$ implies the coexistence scenario.
The frequences and the bifurcation values of the parameter $\tau = \tau^*-\delta_n$ can be approximated by the formulas
\begin{equation}
\omega_n = \frac{2 \pi n}{T} \left(1 - \frac 1 {G_* T \kappa \gamma}\right),
\end{equation}
\begin{equation}
\delta_n =\left(\frac{2 \pi n}{T \gamma}\right)^2 \frac{m \left(\left(\frac{2 \pi n r}{T}\right)^2\frac{1}{\frac{g_0}{G_*} r+G_* m
   \nu}+\frac{g_0 r}{G_*}+G_* m \nu \right)}{2 \kappa
   \left(g_0-\alpha  G_*\right) \left(\left(\frac{2 \pi  n}{T}\right)^2 \frac{r (\mu
   s-f r)}{F_*}-\frac{g_0 r}{G_*}-G_* m \nu \right)}.
\end{equation}
Figure \ref{fig:fast:compete} shows the pulsating behavior of the system \eqref{A1}--\eqref{G1}.

\section{Scaling with $\gamma$. Approximate solution}\label{scaling}
In order to analyze and approximate the asymptotic behavior of the pulsating periodic solution for large $\gamma$, we adapt
the approach proposed by New and Haus for modeling optical systems in \cite{New, Haus} by partial differential equations and an extension of this approach to delay differential models
of mode-locked semiconductor lasers developed in \cite{VladimirovTuraev}.

Consider a pulsating periodic solution of Eqs. \eqref{A}--\eqref{G}.
We divide the period into two stages, the short fast stage $t_b\le t\le t_e$ containing the pulse, when $A$ is large, and the slow stage $t_e \le t\le
t_b+\tau$, during which $A$ is close to zero. Here $\tau\approx T$ is the period of the solution, $t_b$ is the moment when a pulse begins,
$t_e$ is the moment when the pulse ends, $t_e-t_b\ll 1$.

During the pulse of $A$, the terms $AQ$ and $AG$ in the $Q$ and $G$ equations are large compared to the other terms, which therefore can be neglected.
Hence, during the fast stage, Eqs. \eqref{Q}, \eqref{G} can be approximated by the equations
\begin{equation*}
\begin{array}{rcl}
\gamma^{-1}Q'&=&- sAQ,\\
G'&=& -kAG.
\end{array}
\end{equation*}
Integrating, we obtain
\begin{equation}\label{qg}
Q(t)=Q_b e^{-\gamma s P(t)},\qquad G(t)=G_b e^{-k P(t)}
\end{equation}
where $Q_b=Q(t_b)$, $G_b=G(t_b)$ and
\begin{equation}\label{int}
P(t)=\int_{t_b}^t A(\theta)\,d\theta.
\end{equation}
In particular, for the value $G(t_e)=G_e$ at the moment $t=t_e$, we have
\begin{equation}\label{geb}
G_e=G_b e^{-k p}
\end{equation}
where
\begin{equation*}
p=\int_{t_b}^{t_e} A(\theta)\,d\theta.
\end{equation*}
Integrating \eqref{A} over the fast stage and using the fact that $A$ is close to zero at the moments $t_b$ and $t_e$, we obtain the approximate equation
\begin{equation}\label{p}
p=\kappa \int_{t_b-T}^{t_e-T} G(\theta) A(\theta)\,d\theta - \mu \int_{t_b}^{t_e} Q(\theta) A(\theta)\,d\theta.
\end{equation}
We choose the interval $t_b\le t\le t_e$ to be wide enough
to ensure that the interval $t_b-T\le t\le t_e-T$ contains the pulse of $A$ which occurs
between the moments $t_b-\tau$ and $t_e-\tau$ (although the period $\tau$ is not exactly equal to $T$).
Hence, the integrals in the right hand side of \eqref{p} are essentially integrals over the successive pulses.
Therefore, using the periodicity of the solution and relations \eqref{qg},
we can rewrite \eqref{p} as
\begin{eqnarray*}
\fl p=\kappa \int_{t_b-\tau}^{t_e-\tau} G(\theta) A(\theta)\,d\theta - \mu \int_{t_b}^{t_e} Q(\theta) A(\theta)\,d\theta=\\ \kappa \int_{t_e}^{t_b} \left(G_b e^{-k P(\theta)}-Q_b e^{-\gamma s P(\theta)}\right) A(\theta)\,d\theta .
\end{eqnarray*}
As $\gamma\gg1$, the term  $Q(t)=Q_b e^{-\gamma s P(t)}$ can be neglected. (Relation $Q(t)=Q_b e^{-\gamma s P(t)}$ suggests that $Q(t)$ drops to zero on a time scale which is faster
then the time scale of the fast species when the pulse arrives, stays near the zero during the pulse,
and recovers to the equilibrium value $q_0/\beta$ on the same faster time scale after the pulse passes, see the plot of the $Q$-component Figure \ref{main_trace}. 
Therefore, using \eqref{int},
\begin{equation}\label{pp}
p=\kappa G_b \int_{t_e}^{t_b} e^{-k P(\theta)} A(\theta)\,d\theta = \kappa G_b \int_{0}^{p} e^{-k P(\theta)} \,d P(\theta)=\frac{\kappa G_b }k (1-e^{-k p}).
\end{equation}

\begin{figure}[H]
\begin{center}
\subfloat[\label{fig:pstar:sol}]{\includegraphics*[width=.5\columnwidth]{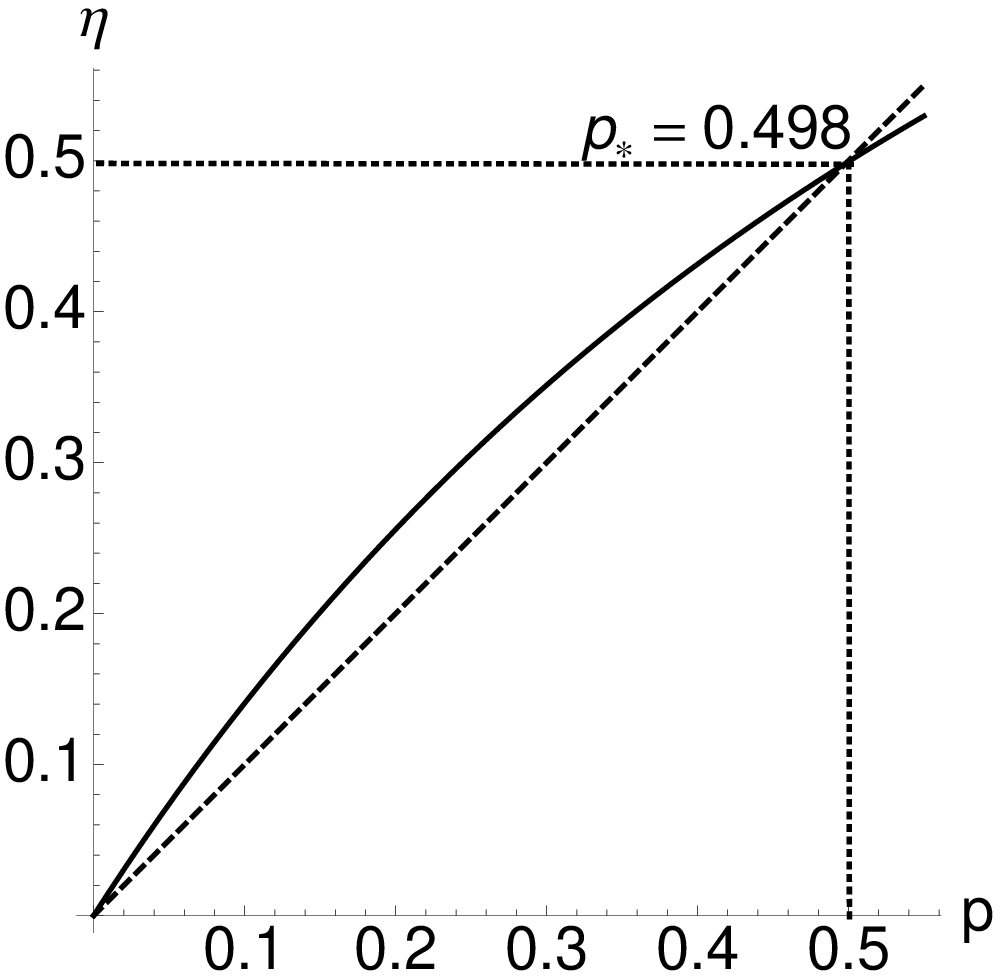}}
\hfill
\subfloat[\label{fig:int:gamma}]{\includegraphics*[width=.5\columnwidth]{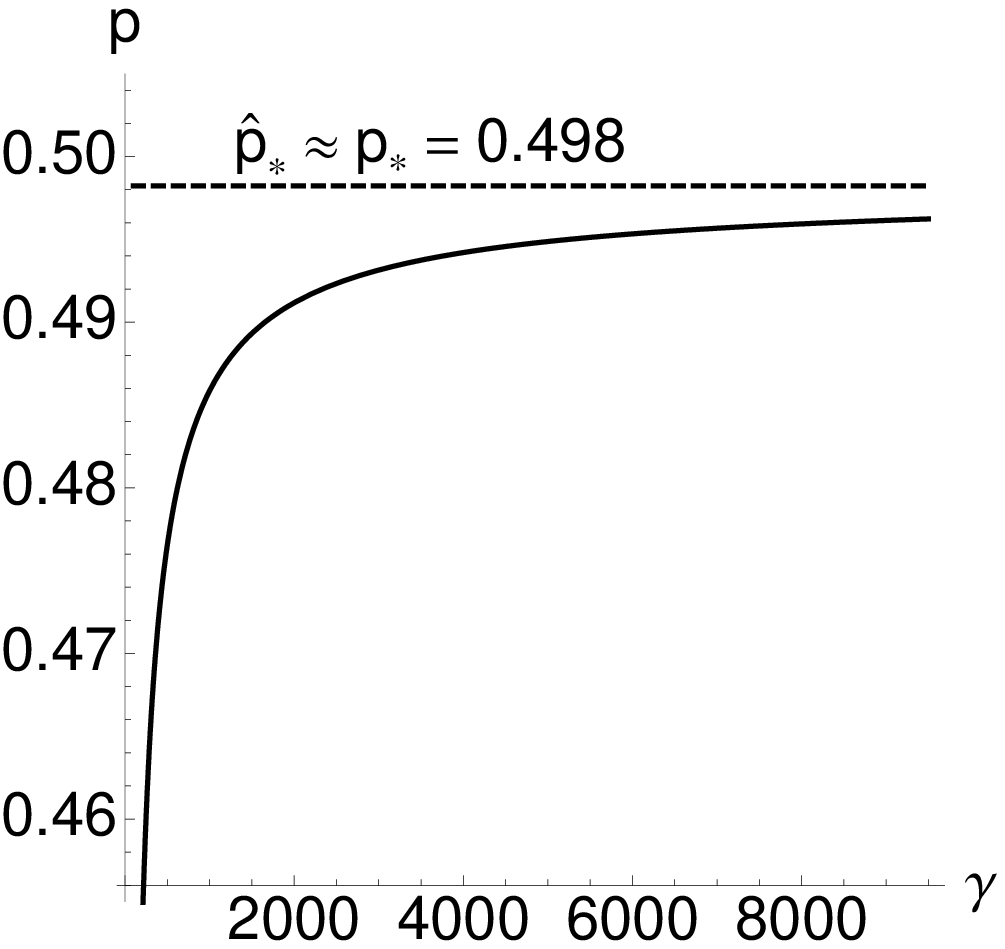}}
\end{center}
\caption{Panel (a) shows the solution of \eqref{p*}. Panel (b) shows the dependence of the integral of $A$-component of system \eqref{A}--\eqref{G} over one period on $\gamma$. The power law fit is shown by the solid line. The horizontal asymptote $p=\hat{p}_*$ coincides with analytic value $p_*$ shown on panel (a).  Here $s = 2$ and $g_0 = 2,6$, all other parameters are the same as in Figure \ref{main_trace}.}
\end{figure}

During the slow stage, the terms $AQ$ and $AG$ are small compared to the other terms in the $Q$ and $G$ equations. Neglecting these terms
results in the linear equations
\begin{equation*}
\begin{array}{rcl}
\gamma^{-1}Q'&=&q_0-\beta Q,\\
G'&=&g_0-\alpha G.
\end{array}
\end{equation*}
Integrating the $G$-equation over the slow stage and combining the integral
\begin{equation*}
g_0- \alpha G_b=(g_0-\alpha G_e)e^{-\alpha(t_b+\tau-t_e)}\approx (g_0-\alpha G_e)e^{-\alpha T}
\end{equation*}
with \eqref{geb}, we obtain
\begin{equation}\label{gb}
G_b = \frac{g_0(1-  e^{-\alpha T})}{\alpha(1- e^{-\alpha T-k p} )}.
\end{equation}
Hence, \eqref{pp} implies
\begin{equation}\label{p*}
p=\frac{\kappa g_0 (1-  e^{-\alpha T})(1-e^{-k p})}{\alpha k (1- e^{-\alpha T-k p} )}=:\eta(p).
\end{equation}
The right hand side $\eta(p)$ is an increasing concave function of $p$ on the positive semiaxis,
which is zero at zero, has the derivative $\kappa g_0/\alpha>1$ at zero, and converges to a constant
as $p\to \infty$, see Figure \ref{fig:pstar:sol}. Therefore, \eqref{p*} has a unique positive root $p_*$, which is the limit of the integral of the $A$-component over a period as $\gamma\to\infty$. Figure \ref{fig:int:gamma} compares the value $p(\gamma)$ of this integral with the limit value $p_*$. The integral has been evaluated numerically for 10 values of $\gamma$ from the interval $100\leq\gamma\leq2000$ by direct simulation of equations \eqref{A}--\eqref{G}. The power law fit
\begin{equation*}
	 \phi(\gamma) = \hat{p}_* + b \gamma^{-\beta},
\end{equation*}
was used to obtain the estimate $\hat{p}_*$ of the limit value $p_*$ of the integral. For the parameter set in Figure \ref{fig:int:gamma}, the error between the numerical estimate $\hat{p}_*$ and the analytic value of $p_*=0.498$ obtained from \eqref{p*} satisfies $\left|p_*-\hat{p}_*\right|<10^{-3}$.

We conclude that in the limit of $\gamma$ tending to infinity,
the $A$-component of the periodic solution converges to the periodic sequence of delta functions (Dirac comb),
\begin{equation}\label{Ato}
A\to p_*\sum_{n=-\infty}^{\infty}\delta (t-nT),
\end{equation}
which has the period $T$ equal to the delay. The $Q$-component is equal to the equilibrium value $q_0/\beta$
between the pulses of $A$ and drops to zero during the pulse. 
The $G$-component grows according to the equation $G'=g_0-\alpha G$ from the value $G_e=G_b e^{-k p_*}$
to the value $G_b$ defined by \eqref{gb}, {\em i.e.},
\begin{equation*}
G(t)=\frac{g_0}\alpha-\left(\frac{g_0}\alpha-G_e\right)e^{-\alpha(t-nT)},\qquad nT<t<(n+1)T,
\end{equation*}
between the pulses of $A$, and drops back to the value $G_e$ during the pulse.
The interval $g_0^*<g_0 <g_0^*+\varepsilon$ of the parameter values, over which the cycle born
via Hopf bifurcation on the positive equilibrium transforms to the pulsating solution, collapses to the threshold, that is $\varepsilon\to0$,
as $\gamma$ grows to infinity. 
In other words, the pulsating solution
described by \eqref{Ato}
can be found "immediately" beyond the threshold for large $\gamma$.

The above approximation does not provide information about the fast stage of the solution such as
the profile of the pulse or the deviation of the period from the delay $T$. In order to obtain such information,
one can adapt the approach of Haus and its modifications, see \cite{Haus, VladimirovTuraev}.
We briefly outline a possible approach without going into much detail here.

Formulas \eqref{bifi} and \eqref{omega1} suggest that the period of the periodic solution scales with $\gamma$ as
\begin{equation}\label{per}
\tau=T\left(1+\frac c{\gamma T}\right)+O(\gamma^{-2}).
\end{equation}
Using this asymptotic formula and the periodicity of the solution,
we can rewrite \eqref{A} as
\begin{equation*}
\gamma^{-1} A'(t-c\gamma^{-1}) + A(t-c\gamma^{-1}) + \mu Q(t-c\gamma^{-1}) A(t-c\gamma^{-1}) = \kappa G(t) A(t).
\end{equation*}
Integrating this equation from $t_b$ over a part of the fast stage $t_b\le t\le t_e$, using \eqref{qg},
and taking into account that $A(t_b)\approx 0$ and $Q$ is almost zero when $A$ becomes large, we obtain an approximate equation
\begin{equation*}
\gamma^{-1} A(t-c\gamma^{-1}) + \int_{t_b}^t A(\theta-c\gamma^{-1})\,d\theta =\kappa G_b \int_{t_b}^t e^{-k P(\theta)} A(\theta)\,d\theta.
\end{equation*}
As $A=P'$, we obtain, similarly to \eqref{pp},
\begin{equation}\label{pulse}
\gamma^{-1} P'(t-c\gamma^{-1})+P(t-c\gamma^{-1})=\frac{\kappa G_b }k (1-e^{-k P(t)}),
\end{equation}
where $G_b$ is defined by \eqref{gb} with $p=p_*$, and $p_*$ is the positive root of \eqref{p*}.
Introducing the fast and reversed time scale $\theta=-\gamma t$, and changing the variable $\bar{P}(\theta)=P(t-c \gamma^{-1})$ we rewrite \eqref{pulse} as
\begin{equation}\label{pulse'}
 -\bar{P}^{\prime}(\theta)+\bar{P}(\theta)=\frac{\kappa G_b }k (1-e^{-k \bar{P}(\theta-c)}),
\end{equation} where $G_b$ is defined by Eqs. \eqref{gb}, \eqref{pp}.
A single pulse of the pulsating periodic solution is, therefore, described by a solution of \eqref{pulse'}
satisfying the boundary conditions
\begin{equation}\label{pulse"}
\bar{P}(-\infty)=0,\qquad \bar{P}(\infty)=p_*.
\end{equation}
Note that both $0$ and $p_*$ are equilibrium points of \eqref{pulse'}.
Therefore conditions \eqref{pulse"} define a heteroclinic orbit of this equation.
%

Linearizing system \eqref{pulse'} {at} 0 and $p_*$ and denoting $a = \kappa G_b$, we obtain
\begin{equation}\label{linp0}
-P^{\prime}(\theta) + P(\theta) = a\, P (\theta-c)
\end{equation} and
\begin{equation}\label{linpstar}
-P^{\prime}(\theta) + P(\theta) = a \, e^{-k p_*} P (\theta-c),
\end{equation} respectively. The characteristic equation of linearization \eqref{linp0} is
\begin{equation}
-\lambda + 1 = a\, e^{-\lambda c} \label{realCharZero}
\end{equation}
or, in real form,
\begin{equation}
\begin{cases}
1 - \alpha  =  a \, e^{-\alpha c} \cos{\beta c}, \\
\beta  =  a \, e^{-\alpha c} \sin{\beta c},
\end{cases}\label{complexCharZero}
\end{equation} where $\lambda = \alpha + i \beta$. Let us assume that
\begin{equation}\label{ac}
a\, c < 1.
\end{equation}
Then the second equation of system \eqref{complexCharZero} implies that $\beta=0$ whenever $\alpha \geq 0$. Furthermore, if $\beta=0$ then the first equation of \eqref{complexCharZero}, condition \eqref{ac} and the relationship ${k p_*}/{(1-e^{-k p_*})}>1$, which follows from \eqref{pp}, imply that $\alpha<0$. Thus, the relation \eqref{ac} guarantees the asymptotic stability of the zero equilibrium of equation \eqref{pulse'}.

The characteristic equation of \eqref{linpstar} has the form
\begin{equation}\label{charLinPStar}
-\lambda + 1 = a \, e^{-k p_*} e^{-\lambda c}.
\end{equation}
Using again the real form of the characteristic equation and the relation $a \, e^{-k p_*} < 1$, which follows from \eqref{pp}, we obtain that \eqref{charLinPStar} has one positive real root and all the other eigenvalues satisfy $\mathrm{Re} \lambda < 0$.
Therefore, the equilibrium $p_*$ of \eqref{pulse'} has a one-dimensional unstable manifold consisting of two nonequilibrium trajectories (and $p_*$ itself). If one of this trajectories belongs to the basin of attraction of the zero equilibrium then this trajectory is a heteroclinic orbit of \eqref{pulse'} which approximates the profile of a pulse for the periodic pulsating solution of Eqs. \eqref{A}--\eqref{G}.

More precisely, if $\tilde P$ denotes the heteroclinic solution of \eqref{pulse'} satisfying \eqref{pulse"},
and $\tilde A=\tilde P'$,
then a pulse of the periodic solution of system \eqref{A}-\eqref{G} is approximated by the formula
\begin{equation*}
A(t)=\gamma \tilde A(\gamma t)
\end{equation*}
for large $\gamma$. Hence, according to this approximation, the amplitude of the pulse scales linearly with $\gamma$, the width of the pulse is inverse proportional to $\gamma$, and the period is approximated by \eqref{per}.

\begin{figure}[ht]
\centering
\includegraphics*[width = 0.6\columnwidth]{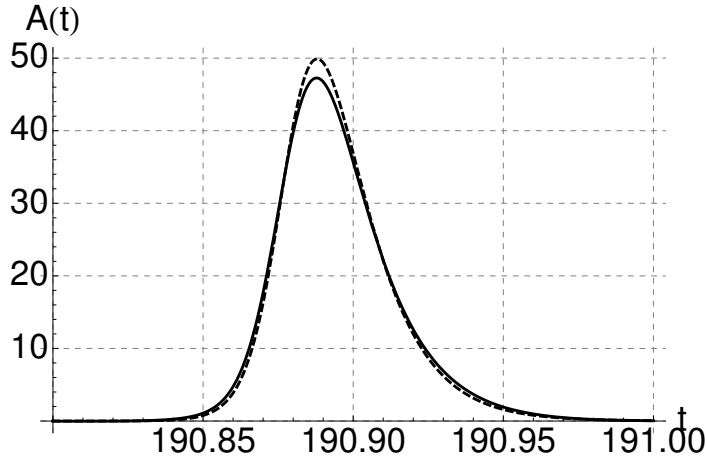}
\caption{Solid line represents one single pulse of the $A$-component of the periodic pulsating solution of system \eqref{A}--\eqref{G} with the following parameters: $g_0 = 5.672, k = 1,q_0 = 1,s = 1.7, T = 1,\alpha = 1, \beta = 1,\gamma = 100,\kappa = 0.6,\mu = 0.5$  (solid). Dashed line is a derivative of the solution of the system \eqref{pulse'} with boundary conditions \eqref{pulse"} where $a = 2.2705, k = 1, c = 0.097$ and $ p_* = 1.9463$.}
\label{fig:pulse1}
\end{figure}

Figure \ref{fig:pulse1} compares this approximation with the numerical solution of system \eqref{A}--\eqref{G}.

In order to solve equation \eqref{pulse'} one needs to know the value of the parameter $c$. If the value of $g_0$ is close the threshold value $g_0^*$ then $c$ can be approximated using equation \eqref{omega1}. For the simulation presented in Figure \ref{fig:pulse1} the value of $c =\gamma (T-\tau)$ was evaluated numerically.

We note that condition \eqref{ac} does not always guarantee that the leading eigenvalue of the zero equilibrium is real as in Figure \ref{fig:pulse1}. If it is complex, then small oscillations of the profile of solution to the approximating equation \eqref{pulse'} are obsered near zero.

The interaction of the three populations in the pulsating regime of Eqs. \eqref{A}--\eqref{G} can be described as follows.
The majority of offsprings of the $A$-species are produced while the population of $A$ is large, that is the pulse of $A$ creates the birth pulse of
their offsprings. This next generation of the $A$-species, however, is not included in the equations, which describe the adults of $A$ only,
until the offspring mature and become part of the system. Thus, the pulse of $A$ creates the next pulse after the delay time $T$.
When this pulse arrives, the $A$ species quickly eliminates the population of the competing $Q$-species. This elimination occurs on a fast time scale. When the population of $Q$ drops to zero, the remaining species $A$ and $G$ exhibit the predator-prey dynamics during the pulse,
{\em i.e.}, over a time interval of the order of $\gamma^{-1}$. The population of the predator $A$ first grows eliminating a large fraction of the prey $G$ and then drops to zero.
Between the pulses, the population of the $A$-species is very close to zero, the population of the $Q$-species quickly returns to its equilibrium value $q_0/\beta$
over a short time interval of the order of $\gamma^{-1}$ in the absence of the competing species $A$, while the population of the $G$-species
slowly recovers between the pulses in the absence of the predator, aiming towards its equilibrium value $g_0/\alpha$, but not having enough time to get close to it.


\section{The role of the competing fast species}\label{competing}
Here we discuss the role of the fast $Q$-species, which competes with the $A$-species,
in creating the pulsating periodic dynamics. This role is critical, even though, as we have seen in the previous section,
important parameters of the periodic pulses are defined by the $A$ and $G$-equations and are independent of the $Q$-equation.
We first discuss the effect of the $Q$-species on the dynamics in terms of the bifurcation scenario considered
in Sections \ref{main}, \ref{variations}, and then in terms of the interaction of the species, along the lines of the analysis presented in Section \ref{scaling}.

In order to highlight the role of the $Q$-species, we compare the dynamics of system \eqref{A}--\eqref{G}
with that of the dynamics of the system
\begin{eqnarray}
\gamma^{-1} A' &=& -A + \kappa G(t-T)A(t-T),\label{2dA}\\
 G' &=& g_0-\alpha G - k A G \label{2dG}
\end{eqnarray}
which is obtained by setting $Q=0$ in \eqref{A} and dropping \eqref{Q}.
Dynamics of system \eqref{2dA}, \eqref{2dG} is essentially the same as dynamics of
system \eqref{A}--\eqref{G} with zero immigration rate $q_0=0$ of the $Q$-species.

System \eqref{2dA}, \eqref{2dG} has two equilibrium points
\begin{equation*}
A=0,\quad G=\frac{g_0}\alpha;\qquad \quad
A=\frac{\kappa (g_0-g_0^*)}k=\frac{\kappa\delta}k,\quad G=\frac1\kappa
\end{equation*}
which collide in the transcritical bifurcation at the threshold value
\begin{equation}\label{thre}
g_0^*=\frac{\alpha}\kappa
\end{equation}
of the bifurcation parameter $g_0$. Like in the case of the three dimensional systems
considered in Sections \ref{main}, \ref{variations},
the equilibrium with zero $A$ is stable below the threshold and unstable above the threshold, while
the equilibrium with nonzero $A$ is positive and stable above the threshold, {\em i.e.}, for $g_0>g_0^*$
(without any additional assumptions about the parameters of Eqs. \eqref{2dA}, \eqref{2dG}).
The unstable equilibrium undergoes the cascade of Hopf bifurcations at the bifurcation points, and with the frequencies,
defined by relations \eqref{transc}, \eqref{transc1} and satisfying the asymptotic formulas \eqref{bifi}, \eqref{bifi:delta}.
However, the positive equilibrium remains stable for all $g_0>g_0^*$ and the system exhibits the equilibrium dynamics rather than a periodic dynamics
above the threshold.
The reason is that the equilibrium with nonzero $A$ undergoes the cascade of Hopf bifurcations below the threshold, that is
in the parameter domain $g_0<g_0^*$ where this equilibrium has a negative $A$-component and is unstable, rather than above the threshold where the equilibrium is positive
and stable. Indeed, substituting the ansatz $\lambda=i\omega$ in the characteristic equation
\begin{equation}\label{jjj}
\gamma^{-1}(\lambda +\kappa g_0) \lambda +(1-e^{-\lambda T})\lambda +\kappa g_0 -\alpha e^{-\lambda T}=0
\end{equation}
of the linearization of the system at the equilibrium with nonzero $A$, we obtain the asymptotic formula
\begin{equation*}
\delta=-\frac{\alpha}{2\kappa}\left(\frac{2\pi n}{\gamma T}\right)^2\left(1+\left(\frac{2\pi n}{\alpha T}\right)^2\right)+O(\gamma^{-3})
\end{equation*}
where the negative sign of $\delta=g_0-g_0^*$ indicates that the Hopf bifurcation occurs below the threshold.
Equation \eqref{jjj} implies
\begin{equation*}
(\omega^2+\kappa^2 g_0^2)(1+\omega^2 \gamma^{-2})=\omega^2+\alpha^2
\end{equation*}
for $\lambda=i\omega$, which is only possible for $g_0\le g_0^*=\alpha/\kappa$, that is below the threshold, thus proving stability
of the positive equilibrium. Clearly \eqref{jjj} cannot have real positive eigenvalues for $g_0>g_0^*$ either.

Note that, in the limit of vanishing $q_0$ where system \eqref{2dA}, \eqref{2dG} approximates system \eqref{A}--\eqref{G},
condition \eqref{trans} is satisfied, while condition \eqref{trans'} is not.

A similar behavior is exhibited by the system of equations \eqref{GG}, \eqref{2dA},
\begin{eqnarray}
\gamma^{-1} A' &=& -A + \kappa G(t-T)A(t-T), \label{compA}\\
 G' &=& g_0 G-\alpha G^2 - k A G, \label{compG}
\end{eqnarray} obtained by setting $Q=0$ in system \eqref{AA}--\eqref{GG} and approximating this system for $q_0=0$.
The equilibrium
\begin{equation}\label{eq+}
A=\frac{g_0-g_0^*}k=\frac\delta{k},\qquad G=\frac1\kappa
\end{equation}
of system \eqref{compA}, \eqref{compG} is positive above the threshold \eqref{thre} and stable on an interval $g_0^*<g_0<g^*_0+\Delta$
of length $\Delta$, which does not vanish for large $\gamma$. Equilibrium \eqref{eq+} is unstable below the threshold, where its $A$-component is negative.
It undergoes the cascade of Hopf bifurcations at the points
\begin{equation*}
\delta=-\frac{\alpha}{2\kappa}\left(\frac{2\pi n}{\gamma T}\right)^2\left(1+\left(\frac{2\pi \kappa n}{\alpha T}\right)^2\right)+O(\gamma^{-3})
\end{equation*}
below the threshold. This equilibrium collides in the transcritical bifurcation, and exchanges stability, with the equilibrium $A=0$, $G=g_0/\alpha$, which
is stable below and unstable above the threshold. Therefore, system \eqref{compA}, \eqref{compG} exhibits starionary dynamics
for all the values $g_0^*<g_0<g_0^*+\Delta$
(this system also has the zero equilibrium, which is always unstable). In this example, the interval of stability of the positive equilibrium is however finite.




%

\section{Conclusion}
We have explored bifurcation scenarios associated with the formation of a pulsating periodic regime in slow-fast
delay differential systems. This regime has a period close to the delay and is characterized by a specific scaling of the pulse width and hight
with the parameter $\gamma\gg 1$ measuring the ratio of the fast and slow time scales. It is formed close to the threshold value
of the bifurcation parameter, at which the zero equilibrium undergoes the transcritical bifurcation and a positive equilibrium appears. Through a case study of
several population dynamics models, we have shown that the formation of the periodic pulsating solution
is associated with a cascade of multiple, almost simultaneous resonant Hopf bifurcations that occur in the immediate vicinity
of the threshold on the positive equilibrium. Using the asymptotic analysis at zero, we have obtained explicit
relationships between the parameters, which ensure this scenario for population dynamics examples.
In particular, we have highlighted the role of competition and shown that the pulses
with the associated cascade of Hopf bifurcations appear when fast species compete
and the interspecific competition is stronger than the intraspecific competition;
in the absence of competition, pulses do not form near the threshold.

We have adapted methods used in laser modeling for separating
 and matching slow and fast stages of the dynamics to obtain an approximation
to the pulsating solution. A modified approach due to New provides
an accurate prediction of the area of the pulse. Furthermore, a modification of the method of Haus has allowed us
to obtain asymptotics of the pulse shape as $\gamma\to\infty$. This shape is described by a heteroclinic solution
of a scalar delay equation that depends only on two parameters. The positive equilibrium of
this equation has a one-dimensional unstable manifold, which, as numerical simulations show, is attracted to the stable
zero equilibrium, thus forming the hetroclinic connection that describes the profile of the pulse.

Population dynamics setting differs from that of laser dynamics in several respects.
In particular, population models are positively invariant and the pulsating regime in this setting is positive; typically, these models contain
quadratic nonlinearities as opposed to cubic (or exponential) nonlinearities of laser models;
finally, the pulsating variable, which plays the role of the counterpart of $A$ in laser dynamics, is complex-valued.

Due to positive invariance, the transcritical bifurcation with the associated zero eigenvalue
is an important ingredient of the bifurcation scenario described in this work.
This scenario can be compared to the Eckhaus and modulational instabilities, which
are well known in the context of spatially distributed systems and have been recently
studied for systems containing large delays \cite{eckhauswy, modinstab}.
The evolution of the pseudocontinuous spectrum of the zero equilibrium shown in Fugure \ref{stabil}a
is similar to the picture associated with the Eckhaus instability. The ``parabola'' carrying the
pseudocontinuous spectrum moves as a whole to the unstable half-plane as the bifurcation parameter increases. Furthermore, as in the Eckhaus scenario \cite{eckhaus}, we observe
the appearance of multiple unstable periodic solutions, which then stabilize
via secondary bifurcations leading to co-existence of multiple periodic attractors, see Figure \ref{figDiag}.
On the other hand, the evolution of the spectrum of the positive equilibrium
that intersects the zero equilibrium in the transcritical bifurcation
reminds the modulation instability scenario, in which the ``parabola'' carrying the pseudocontinuous spectrum develops two humps
that cross the imaginary axis, while the vertex of the parabola at zero is not moving \cite{modinstab}.
Interestingly, although similar humps are observed in Figure \ref{stabil1}, they are formed through a different mechanism.
Namely, eigenvalues with smaller imaginary part that belong to the pseudocontinuous spectrum get absorbed by the strongly stable spectrum
as the bifurcation parameter increases. This interaction of the pseudocontinuous and strongly stable spectra
results in the formation of humps and, further, in stabilization of the positive equilibrium for higher values
of the bifurcation parameter. However, 
a common feature of all the above scenarios
 is that eigenvalues with smaller imaginary part cross the imaginary axis from the stable to the unstable domain
before eigenvalues with larger imaginary part do. Hence, all these scenarios can be viewed 
as longwavelength instabilities.


\section*{Acknowledgments}
P. K. and D. R. acknowledge the support of NSF through grant DMS-1413223. A. V. was partially supported by grant RSF 14-41-00044.

\section*{References}
\bibliographystyle{unsrt}
\bibliography{biblio}

\end{document}